\documentclass{amsart}
\usepackage{amscd}
\numberwithin{equation}{section}
\let\cal\mathcal

\newtheorem{lemma}{Lemma}[section]
\newtheorem{proposition}[lemma]{Proposition}
\newtheorem{theorem}[lemma]{Theorem}
\newtheorem{theoremA}[lemma]{Theorem A}
\newtheorem{theoremB}[lemma]{Theorem B}

\newtheorem{corollary}[lemma]{Corollary}

\newtheorem{fact}[lemma]{Fact}

\theoremstyle{definition}

\newtheorem{example}[lemma]{Example}
\newtheorem{definition}[lemma]{Definition}

\newtheorem{notation}[lemma]{Notation}
\newtheorem{conventions}[lemma]{Conventions}

{

\newtheorem{step}{Step}

\newtheorem*{proofoftheorem}{Proof of  Theorem 2.20}
\newtheorem*{prooftheorem}{Proof of Theorem 3.1}
\newtheorem*{prooftheoremA}{Proof of Theorem A}
\newtheorem*{prooftheoremB}{Proof of Theorem B}

\newtheorem*{proofproposition}{Proof of Proposition 2.21}
\newtheorem*{proofproposition3}{Proof of Proposition 3.10}
\newtheorem*{prooflemma}{Proof of lemma}
}

\theoremstyle{remark}

\newtheorem{remark}[lemma]{Remark}

\newcommand{\Acal}{\mbox{$\cal A$}}
\newcommand{\Bcal}{\mbox{$\cal B$}}

\newcommand{\Dcal}{\mbox{$\cal D$}}

\newcommand{\Fcal}{\mbox{$\cal F$}}

\newcommand{\Lcal}{\mbox{$\cal L$}}

\newcommand{\Ncal}{\mbox{$\cal N$}}

\newcommand{\Ocal}{\mbox{$\cal O$}}

\newcommand{\Wcal}{\mbox{$\cal W$}}
\newcommand{\Rcal}{\mbox{$\cal R$}}

\newcommand{\Scal}{\mbox{$\cal S$}}

\newdimen\uboxsep \uboxsep=1ex
\def\uboxn#1{\vtop to 0pt{\hrule height 0pt depth 0pt\vskip\uboxsep
\hbox to 0pt{\hss #1\hss}\vss}}

\def\uboxs#1{\vbox to 0pt{\vss\hbox to 0pt{\hss #1\hss}
\vskip\uboxsep\hrule height 0pt depth 0pt}}

\title[Binomial solutions of the Yang-Baxter equation]{Binomial skew polynomial rings, Artin-Schelter regularity, and
binomial solutions of the Yang-Baxter equation} 
\thanks{The author was partially supported by the Department of Mathematics of  Harvard University,
by Grant MM1106/2001 of the Bulgarian National Science Fund
  of the Ministry of Education and Science,
  and by The Abdus Salam International Centre
   for Theoretical Physics (ICTP)}
\author{Tatiana Gateva-Ivanova}
\address{Institute of Mathematics and Informatics\\
Bulgarian Academy of Sciences\\
Sofia 1113, Bulgaria
and\\
American University in Bulgaria\\
2700 Blagoevgrad, Bulgaria }

\email{tatianagateva@yahoo.com, tatyana@aubg.bg}

\begin{document}
%\date{\today}
\begin{abstract}
Let $\textbf{k}$ be a field and $X$ be a set of $n$ elements. We introduce
and study a class of quadratic $\textbf{k}$-algebras
called \emph{quantum binomial algebras}. Our main result  shows
that such an
algebra $A$ defines a solution of the classical Yang-Baxter
equation  (YBE), if and only if its Koszul dual $A^{!}$ is
Frobenius of dimension $n,$ with a \emph{regular socle} and for
each $x,y \in X $ an equality of the type $xyy=\alpha zzt,$ where
$\alpha \in k \setminus\{0\},$ and $z,t \in X$ is satisfied in
$A$. We prove the equivalence of the notions \emph{a binomial skew
polynomial ring} and \emph{a binomial solution of YBE}. This
implies that the  Yang-Baxter algebra of such a solution is of
Poincar\'{e}-Birkhoff-Witt type, and possesses a number of other
nice properties such as being Koszul, Noetherian, and an
Artin-Schelter regular domain.
\end{abstract}
\maketitle
\section{Introduction}

   In the paper we work with associative finitely presented
graded $\textbf{k}$-algebras $A=\bigoplus_{i\ge0}A_i$, where
$\textbf{k}$ is a field,  $A_0=\textbf{k}$, and $A$ is generated
by $A_1.$
%therefore each $A_i$ is finite dimensional.
We restrict
our attention to  a class of algebras with quadratic binomial
defining relations and study the close relations between different
algebraic notions such as \emph{Artin-Schelter regular rings},
\emph{Yang-Baxter algebras} defined via \emph{binomial solutions} of the
classical Yang-Baxter equation, and a class of quadratic standard
finitely presented algebras with a Poincar\'{e}-Birkhoff-Witt type
$\textbf{k}$-basis, called \emph{binomial skew polynomial rings}.

Following a classical tradition (and recent trend), we take a
combinatorial approach to study $A$. The properties of $A$ will be
read off a presentation $A= \textbf{k} \langle X\rangle /(\Re)$,
where $X$ is a finite set of indeterminates of degree $1$,
$\textbf{k}\langle X\rangle$ is the unitary free associative
algebra generated by $X$, and $(\Re)$ is the two-sided ideal of
relations, generated by a {\em finite} set $\Re$ of homogeneous
polynomials.

Artin and Schelter \cite{AS} call a graded algebra $A$ \emph{ regular} if
\begin{enumerate}
\item[(i)]
$A$ has {\em finite global dimension\/} $d$, that is, each graded
$A$-module has a free resolution of length at most $d$.
\item[(ii)]
$A$ has \emph {finite Gelfand-Kirillov dimension}, meaning that
the integer-valued function $i\mapsto\dim_{\textbf{k}}A_i$ is
bounded by a polynomial in $i$.
\item[(iii)]
$A$ is \emph{Gorenstein}, that is, $Ext^i_A(\textbf{k},A)=0$ for
$i\ne d$ and $Ext^d_A(\textbf{k},A)\cong \textbf{k}$.
\end{enumerate}

The regular rings were introduced  and studied first in \cite{AS}.
When $d\le3$ all regular algebras are classified. The problem of
classification of regular rings is difficult and remains open even
for regular rings of global dimension $4$. The study of
Artin-Schelter regular rings, their classification, and finding
new classes of such rings is one of the basic problems for
noncommutative geometry. Numerous works on this topic appeared
during the last 16 years, see for example \cite{ATM}, \cite{BSM},
\cite{Levasseur}, \cite{paul}, \cite{Michel-Tate},
\cite{Michel-Mart}, etc.

For the rest of the paper we fix $X$. If an enumeration $X =
\{x_{1},\dots, x_{n}\}$ is given, we will consider the
degree-lexicographic order $\prec$ on $\langle X \rangle$, the
unitary free semigroup generated by $X$ (we  assume $x_1 \prec x_2
\prec \cdots \prec x_n$).

Suppose the algebra $A$ is given with a finite presentation  $A =
\textbf{k} \langle x_{1},\dots,x_{n}\rangle /(\Re).$

In some cases we will ignore a given enumeration on $X$ and will
search for an appropriate enumeration (if any), which provides a
degree-lexicographic ordering $\prec$ with respect to which the
relations $\Re$ become of \emph{skew-polynomial type}, see Definition
\ref{binomialringdef}.

Recall that a monomial $u\in \langle X \rangle$ is {\em normal} mod $\Re$
(with respect to the chosen order) if $u$ does not contain as a
segment any of the highest monomials of the polynomials in $\Re$.
By $N(\Re)$ we denote the
set of all normal mod $\Re$ monomials.

\begin{notation} As usual, we denote
$\textbf{k}^{\times}=\textbf{k}\setminus\{0\}$. If  $\omega = x_{i_1}
\cdots x_{i_m} \in \langle X \rangle,$ its length $m$ is denoted
by  $\mid \omega \mid $. $X^m$ will denote the set of all words of
length $m$ in the free semigroup $\langle X \rangle.$ We shall
identify  the $m$-th tensor power $V^{\otimes m}$ with $V^m=
Span_\textbf{k}{X^m} $,  the $\textbf{k}$-vector space spanned by
all monomials of length $m$.
\end{notation}
We shall introduce now a  class of quadratic algebras  with
binomial relations, we call them \emph{quantum binomial algebras}, which
contains various algebras, such as binomial skew polynomial rings,
\cite{T1}, \cite{T96}, \cite{T3}, the Yang-Baxter algebras defined
via the so called \emph{binomial solutions} of the classical
Yang-Baxter equation, \cite{T06}, the semigroup algebras of
semigroups of skew type, \cite{TJO}, etc. all of which are
actively studied.

\begin{definition}
\label{quantumquadraticdef} Let $A(\textbf{k}, X, \Re) =
\textbf{k}\langle X \rangle/(\Re ) $ be a finitely presented
{\textbf{k}}-algebra with a set of generators $X$ consisting of
$n$ elements, and quadratic defining relations $\Re\subset
\textbf{k}\langle X \rangle$. The relations $\Re$ are called
\emph{quantum binomial relations} and  $A$ is \emph{a quantum
binomial algebra} if the following conditions hold.
\begin{enumerate}
\item[(a)]
\label{quantumquadraticdefa} Each relation in  $\Re$ is of the
shape $ xy-c_{yx}y^{\prime}x^{\prime},$ where $x, y, x^{\prime},
y^{\prime} \in X,$ and $c_{xy} \in \textbf{k}^{\times}$ (this is what we
call \emph{a binomial relation}).
\item[(b)]
\label{quantumquadraticdefb} Each $xy, x \neq y$ of length $2$
occurs at most once in $\Re$.
\item[(c)]
\label{quantumquadraticdefc} Each relation is \emph{square-free},
i.e. it does not contain a monomial of the shape $xx,$ $x \in X.$
\item[(d)]
\label{quantumquadraticdefd} The relations $\Re$ are \emph{non
degenerate}, i.e. the canonical bijection $r=r(\Re): X^2
\longrightarrow X^2,$ associated with $\Re$, see Definition
\ref{quantumquadraticmap}  is left and right non degenerate.
\end{enumerate}
A quantum binomial algebra $A$ is called \emph{standard quantum
binomial algebra} if the set $\Re$ is a Gr\"obner basis with
respect to the degree-lexicographic ordering $\prec$,  where  some
appropriate enumeration of $X$ is chosen, $X = \{x_1 \prec x_2
\prec \cdots \prec x_n \}.$
\end{definition}
\begin{definition}
\label{quantumquadraticmap} Let  $\Re\subset \textbf{k}\langle X
\rangle$ be a set of quadratic binomial relations, satisfying
conditions (a)  and (b) of Definition \ref{quantumquadraticdef}.  \emph{The
automorphism associated with} $\Re,$ $R=R(\Re): V^2
\longrightarrow V^2,$ is defined as follows: on monomials which
occur in some relation, $xy-c_{xy}y^{\prime}x^{\prime} \in \Re,$
we set
\[
R(xy)=c_{xy}y^{\prime}x^{\prime}, \quad \text{ and}\quad
R(y^{\prime}x^{\prime})= (c_{xy})^{-1}xy.\]

If $xy$, does not occur in any relation ($x=y$ is also possible),
then we set $R(xy)= xy.$

We also define a bijection $r= r(\Re): X^2 \longrightarrow X^2$ as
$r(xy)=y^{\prime}x^{\prime},$ and $r(y^{\prime}x^{\prime})= xy,$
if $xy-c_{xy}y^{\prime}x^{\prime} \in \Re.$ If $xy$, does not
occur in any relation then we set $r(xy)= xy.$ We call $r(\Re)$
\emph{the (set-theoretic) canonical map associated with} $\Re$.

We say that $r$ is \emph{nondegenerate}, if the two maps $\Lcal_x:
X\longrightarrow X,$ and $\Rcal _y: X\longrightarrow X$ determined
via the formula:
\[
r(xy)=\Lcal _x(y)\Rcal _y(x)
\]
are bijections for each $x, y \in X.$

$R$ is called \emph{non-degenerate} if $r$ is non-degenerate. In
this case we shall also say that the defining relations $\Re$ are
\emph{non degenerate binomial relations}.
\end{definition}
\begin{definition}
\label{associatedsemigroup} With each quantum binomial set of
relations $\Re$ we associate a set of semigroup relations $\Re
_0,$ obtained by setting $c_{xy}=1,$ for each relation
$(xy-c_{yx}y^{\prime}x^{\prime})\in \Re$. In other words,
\[
 \Re _0 =\{ xy=y^{\prime}x^{\prime}  \mid
xy-c_{xy}y^{\prime}x^{\prime} \in \Re \}
\]
The \emph{semigroup associated to} $A(\textbf{k}, X, \Re)$ is
$\Scal_0=\Scal_0(X, \Re_0)=
\langle X ; \Re_0\rangle $, we also refer to it as \emph{quantum
binomial  semigroup}. The semigroup algebra associated to
$A(\textbf{k}, X, \Re)$ is $\Acal _0 = \textbf{k}\langle X
\rangle/(\Re _0)$, which is isomorphic to $ \textbf{k}\Scal _0.$
\end{definition}
The following lemma gives more precise description of the
relations in a quantum binomial algebra. We give the proof in
Section 2.
\begin{lemma}
\label{binomnchoose2relations} Let $A(\textbf{k}, X, \Re)$ be a
quantum binomial algebra, let $\Scal _0$ be the associated quantum
binomial semigroup. Then
the following conditions are satisfied.
\begin{enumerate}
\item[(i)]
$\Re$ contains precisely $\binom{n}{2}$ relations
\item[(ii)]
Each monomial $xy \in X^2, x \neq y,$ occurs exactly once in
$\Re$.
\item[(iii)]
$xy -c_{yx} y^{\prime}x^{\prime}\in \Re,$ implies $y^{\prime} \neq
x,  x^{\prime}\neq y.$
\item[(iv)]
The left and right Ore conditions, (see Definition \ref{Oreconditions}) are
satisfied in $\Scal _0.$
\end{enumerate}
\end{lemma}
Clearly, if the set $\Re$ is a Gr\"{o}bner basis then $\Re _0$ is
also a Gr\"{o}bner basis. Therefore,  for a standard quantum
binomial algebra  $A(\textbf{k}, X, \Re)$ the associated semigroup
algebra $\Acal _0$ is also standard quantum binomial.

\begin{example}
a) Each binomial skew polynomial ring, see Definition \ref{binomialringdef} is a
standard quantum binomial algebra.

b) Let $R$ be a binomial solution of the classical Yang-Baxter
equation, see Definition \ref{binomialsolution}, and let $\Re(R)$  be the
corresponding quadratic relations, then the Yang-Baxter algebra
$A= \textbf{k}\langle X \rangle/(\Re )$ is a quantum binomial
algebra.

c) $A= k\langle x_1, x_2, x_3, x_4\rangle/(x_4x_3- x_2x_4,
x_4x_2-x_1x_3, x_4x_1-x_3x_4, x_3x_2-x_2x_3, x_3x_1-x_1x_4, x_2x_1
- x_1x_2)$ is a quantum binomial algebra, which is not standard
quantum binomial, i.e. whatever enumeration on $X$ we fix, the set
of relations $\Re$ is not a Gr\"{o}bner basis with respect to
$\prec$. This can be deduced by direct computations, but one needs
to check all possible, $4!$ enumerations of $X$, which is too
long. (In particular if we chose the given enumeration, the
ambiguity $x_4x_3x_1$ is not solvable). Here we give another
proof, which is universal and does not depend on the enumeration.
Assume, on the contrary, $\Re$ is a Gr\"{o}bner basis, with
respect to an appropriate enumeration. Therefore $A$ is a binomial
skew polynomial ring and the
%(weak)
cyclic condition is satisfied,
see Definition \ref{weakcycliccondition}. Now the relations $x_4x_3-
x_2x_4,x_4x_2-x_1x_3,$ give a contradiction.
\end{example}

We single out an important subclass of standard quantum binomial
algebras with a  Poincar\'{e}-Birkhoff-Witt type
$\textbf{k}$-basis, namely  \emph{the binomial skew polynomial
rings}. These rings were introduced and studied in \cite{T1},
\cite{T96}, \cite{T3}, \cite{TM}, \cite{Laffaille}. Laffaille calls
them \emph{quantum binomial algebras}.  He shows in
\cite{Laffaille}, that for $\mid X \mid \leq 6$,  the associated
automorphism $R$ is a solution of the Yang-Baxter equation.   We
prefer to keep the name "binomial skew polynomial rings" since we
have been using this name for already 10 years. It was proven in
1995, see \cite{T3} and \cite{TM} that the binomial skew
polynomial rings provide a new  (at that time) class of
Artin-Schelter regular rings of global dimension $n$, where $n$ is
the number of generators $X.$ We recall now the definition.
\begin{definition}\label{binomialringdef}
\cite{T96} A {\em binomial skew polynomial ring\/} is a graded
algebra $A=\textbf{k} \langle x_1, \cdots , x_n\rangle/(\Re)$ in
which the indeterminates $x_i$ have degree $1$, and which has
precisely $\binom{n}{2}$ defining relations $\Re=\{x_{j}x_{i} -
c_{ij}x_{i^\prime}x_{j^\prime}\}_{1\leq i<j\leq n}$ such that
\begin{enumerate}
\item[(a)] \label{binomialringa} $c_{ij} \in \textbf{k}^{\times}$;
%=\textbf{k}\setminus\{0\}$;
\item[(b)]
\label{binomialringb} For every pair $i, j $ $1\leq i<j\leq n$,
the relation $x_{j}x_{i} - c_{ij}x_{i'}x_{j'}\in \Re,$ satisfies
$j
> i^{\prime}$, $i'\leq j'$;
\item[(c)] \label{binomialringc} Every ordered monomial $x_ix_j,$
with $1 \leq i < j \leq n$ occurs in  the right hand side of some
relation in $\Re$;
\item[(d)] \label{binomialringd} $\Re$ is the
{\it reduced Gr\"obner basis\/}
of the two-sided ideal $(\Re)$,(with respect to the order $\prec$ on
$\langle X \rangle$) or equivallently the ambiguities $x_kx_jx_i,$
with $k>j>i$ do not give rise to new relations in $A.$
\end{enumerate}

We call $\Re$ \emph{relations of skew-polynomial type} if
conditions \ref{binomialringdef} (a), (b) and (c) are satisfied
(we do not assume (d)) .
\end{definition}
By \cite{B} condition \ref{binomialringdef} (d) may be rephrased
by saying that  \emph{the set of ordered monomials}
\[
\label{ncal} \Ncal _0 = \{x_{1}^{\alpha_{1}}\cdots
x_{n}^{\alpha_{n}}\mid \alpha_{n} \geq 0 \text{ for } 1 \leq i
\leq n\}
\]
is a $\textbf{k}$-basis of $A$.
\begin{definition}
\label{skewpolynomialsemigroup}
We say that the semigroup $\Scal _0$ is \emph{a
semigroup of skew-polynomial type}, (or shortly, \emph{a
skew-polynomial semigroup}) if it has a standard finite
presentation as $\Scal _0 = \langle X; \Re_0 \rangle $, where the
set of generators $X$ is ordered: $x_1 \prec x_2 \prec \cdots
\prec x_n,$ and the set
\[
\Re _0 =\{ x_jx_i=x_{i^{\prime}}x_{j^{\prime}})\mid 1 \leq i< j
\leq n, 1 \leq i^{\prime} < j^{\prime} \leq n \},
\]
contains precisely $\binom{n}{2}$ quadratic square-free binomial
defining relations, each of them satisfying the following
conditions:
\begin{enumerate}
\item[(i)] \label{skewpolynomialsemigroupi} each monomial $xy\in
X^2$, with $x\neq y$, occurs in exactly one relation in $\Re _0$;
a monomial of the type $xx$ does not occur in any relation in $\Re
_0$; \item[(ii)] \label{skewpolynomialsemigroupii} if
$(x_jx_i=x_{i^{\prime}}x_{j^{\prime}})\in \Re _0$, with $1 \leq i<
j \leq n,$ then  $i^{\prime} < j^{\prime}$, and $j > i^{\prime}$.

(further studies show that this also implies $i < j^{\prime}$
see \cite{T96})
\item[(iii)]
\label{skewpolynomialsemigroupiii}
the monomials $x_kx_jx_i$ with $k>j>i, 1\leq i,j,k, \leq n$ do
not give rise to new relations in $S_0$, or equivalently,  cf.
\cite{B}, $\Re _0$ is a Gr\"{o}bner basis with respect to the
degree-lexicographic ordering of the free semigroup $\langle X
\rangle$.
\end{enumerate}
\end{definition}
\begin{example}
\[A_1 = \textbf{k}\langle x_{1}, x_{2}, x_{3} \rangle /(\Re_1),\]
where
\[\Re_1 = \{x_{3}x_{2} - x_{1}x_{3},\ x_{3}x_{1} - x_{1}x_{3}, \
x_{2}x_{1} - x_{1}x_{2} \}.\] Then $\Re _1$ is a Gr\"obner basis,
but it does not satisfy  (c) in Definition \ref{binomialringdef}, hence $A_1$ is not a binomial skew polynomial ring.
Respectively,   the semigroup $\langle X\mid \Re _0\rangle $
is not a
skew-polynomial semigroup. (Here $\Re _0$ are the associated semigroup relations as in Definition \ref{associatedsemigroup}.
\end{example}
\begin{example} Let
\[ A_2 = \textbf{k}\langle x_{1}, x_{2}, x_{3}, x_{4} \rangle
/(\Re_2),
\]
where
\begin{align}
\Re _2 = \{& x_{4}x_{3} - a x_{3}x_{4},\ x_{4}x_{2} -
bx_{1}x_{3},\ x_{4}x_{1} - cx_{2}x_{3},\notag
\\
& x_{3}x_{2} - dx_{1}x_{4},\ x_{3}x_{1} - ex_{2}x_{4},\ x_{2}x_{1}
- fx_{1}x_{2}\},\notag
\end{align}
and the coefficients $a,b,c,d,e,f$ satisfy
\[
abcdef \not= 0,\ a^{2} = f^{2} = be/cd = cd/be, a^4=f^4=1.
\]
This is a binomial skew polynomial ring.  $A_2$ is regular and
left and right Noetherian domain.
\end{example}
A classification  of the binomial skew polynomial rings with $4$
generators was given in \cite{T1}, some of those algebras are
isomorphic. A computer programme was used in \cite{Laffaille} to
find all the families of binomial skew polynomial rings in the
case $n \leq6$, some of the algebras there  are also isomprphic.
One can also find a classification of the binomial skew polynomial
rings with $5$ generators and various examples of such rings in
$6$ generators found independently in \cite{Coqui}.

Now we recall the definition of the Yang-Baxter equation.

Let $V$ be a vector space over a field $\textbf{k}$. A linear
automorphism $R$ of $V\otimes V$ is  \emph{a solution of the
Yang-Baxter equation}, (YBE)  if the equality
\begin{equation}
\label{YBE} R^{12}R^{23}R^{12} = R^{23}R^{12}R^{23}
\end{equation}
%(R\otimes id_V)(id_V\otimes R)(R\otimes id_V) =
%(id_V\otimes R)(R\otimes id_V)(id_V\otimes R)

holds in the automorphism  group of $V\otimes V\otimes V,$ where
$R^{ij}$ means $R$ acting on the i-th and j-th component.

In 1990 V. Drinfeld \cite{D}  posed the problem of studying  the
\emph{set-theoretic solutions}of YBE.
\begin{definition}
\label{settheoreticYBE} A bijective map $r: X^2 \longrightarrow
X^2$,
is called a \emph{set-theoretic
solution of the Yang-Baxter equation} (YBE) if  the braid relation
\[
\label{settheoreticYBE} r^{12}r^{23}r^{12} = r^{23}r^{12}r^{23}
\]
holds in $X^3,$ where the two bijective maps  $r^{ii+1}: X^3
\longrightarrow X^3$, $ 1 \leq i \leq 2$ are defined as  $r^{12} =
r\times Id_X$, and $r^{23}=Id_X\times r$.

We use notation $(X, r)$ for nondegenerate involutive
set-theoretic solutions of YBE. (For nondegeneracy, see Definition
\ref{quantumquadraticmap}).

Each set-theoretic solution $r$ of the Yang-Baxter equation
induces an operator $R$  on $V\otimes V$ for the vector space $V$
spanned by $X$, which is, clearly, a solution of  \ref{YBE}.
\end{definition}
\begin{definition}
\label{binomialsolution} (\cite{T06}, Def. 9.1) Let $V$ be a finite dimensional vector space over a field
$\textbf{k}$ with a  $\textbf{k}$-basis $X = \{x_1, \cdots, x_n \}$.
Suppose the linear automorphism $R: V\otimes V \longrightarrow V\otimes V$ is a solution of the Yang-Baxter equation.

We say that $R$ is
\emph{a binomial solution of the (classical)  Yang-Baxter equation}
or shortly \emph{binomial solution}  if the following conditions
hold:
\begin{enumerate}
\item for every pair $i\neq j, 1 \leq i,j \leq n,$
\[
 R(x_j\otimes x_i) = c_{ij} x_{i^{\prime}}\otimes
x_{j^{\prime}},
  R(x_{i^{\prime}}\otimes
x_{j^{\prime}}) = \frac{1}{c_{ij}} x_j\otimes x_i,\ \text{where}\  c_{ij}
\in k, c_{ij}\neq 0.
\]
\item $R$ is \emph{non-degenerate}, that is \emph{the associated
set-theoretic solution } $(X, r(R))$, is \emph{non-degenerate},
where $r=r(R): X^2 \longrightarrow X^2$ is defined as
\[
 r(x_jx_i) = x_{i^{\prime}}x_{j^{\prime}}, \;
r(x_{i^{\prime}}x_{j^{\prime}}) = x_jx_i \;
\text{if} \ R(x_j\otimes x_i) = c_{ij} x_{i^{\prime}}\otimes
x_{j^{\prime}},
\]
see \ref{quantumquadraticmap}, see also \cite{ESS}, \cite{T06}.
\item
We call the binomial solution $R$ (respectively the set-theoreric
solution $(X,r)$)  \emph{square-free}   if $R(x\otimes x) =
x\otimes x$, (respectively $r(xx)=xx$) for all $x \in X$
\end{enumerate}
\end{definition}
\begin{notation} By $(\textbf{k}, X, R)$ we shall denote a square-free
binomial solution of the classical Yang-Baxter equation.
\end{notation}
Each binomial solution $(\textbf{k}, X, R)$ defines a
quadratic algebra $\Acal_R= \Acal(\textbf{k}, X, R)$, namely
\emph{the associated Yang-Baxter algebra}, in the sense of Manin
\cite{Maninpreprint}, see also \cite{T06} . The algebra
$\Acal (k, X, R)$ is generated by
$X$ and has quadratic defining relations, $\Re(R)$ determined by
$R$ as in (\ref{defrelations2}):
\begin{equation}
\label{defrelations2} \Re(R) =\{
(x_jx_i-c_{ij}x_{i^{\prime}}x_{j^{\prime}})  \mid R(x_j\otimes
x_i) = c_{ij} x_{i^{\prime}}\otimes x_{j^{\prime}} \}
\end{equation}

Given a set-theoretic solution $(X, r),$ we define the quadratic
relations $\Re(r),$ \emph{the associated Yang-Baxter semigroup}
$S(X, r)$ and the algebra $A(\textbf{k}, X, r)$ analogously, see
\cite{T06}.

\begin{definition} \cite{T06}
\label{weakcycliccondition} Let $A= k\langle X\rangle/(\Re)$ be a
quantum binomial algebra, let $\Scal _0$ be the associated
semigroup. We say that $A,$ respectively $\Scal_0$ satisfies
\emph{the weak cyclic condition} if for any $x, y \in X, x \neq y$
the following relations  hold in $\Scal _0:$
\[
 (yx=x_1y_1)  \in \Re _0 \;
\text{implies}\; (yx_1=x_2y_1) \in \Re _0, \;(y_1x=x_1y_2)  \in
\Re _0 .
\]
for some appropriate $x_2, y_2 \in X.$ Or equivalently, for all
$x, y \in X$ one has<.
\[
\label{weakcyclicconditionrelations} \Rcal_{\Lcal_y(x)}(y)
=\Rcal_x(y), \; \Lcal_{\Rcal_x(y)}(x) =\Lcal_y(x).
\]
\end{definition}
It is shown in \cite{T96} that every binomial skew polynomial ring
$A$  satisfies the weak cyclic condition. Furthermore, every
Yang-Baxter semigroup $S(X,r)$ associated with a set-theoretic
solution $(X,r)$ satisfies the weak cyclic condition,
\cite{TLovetch} and \cite{T06}.
\begin{remark}
\label{cyclicconditionremark}
In fact both $A$ and $S(X,r)$ satisfy a stronger condition which we
call \emph{the cyclic condition}, see \cite{T96}, and \cite{T06}.
\end{remark}
For the main results we need to recall the definitions of the
Koszul dual algebra and of a Frobenius algebra.

The Koszul dual $A^{!}$ is defined in  \cite{Maninpreprint, 5.1}.
One can deduce from there the following presentation of $A^{!}$ in
terms of generators and relations.
\begin{definition}
\label{koszuldualdef} Suppose $A= \textbf{k}\langle X \rangle/(\Re
),$ is a quantum binomial algebra. \emph{The Koszul dual $A^{!}$
of $A$},  \cite{Maninpreprint} is the quadratic algebra,
\[
\textbf{k}\langle \xi _1, \cdots, \xi _n \rangle/(\Re^{\bot}),
\]
where the set $\Re^{\bot}$ contains precisely  $\binom{n}{2}+n$
relations  of the following two types:

a) binomials
\[
\xi_j\xi_i+(c_{ij})^{-1}\xi_{i^{\prime}}\xi_{j^{\prime}}\in
\Re^{\bot}, \ \text{whenever} \
x_jx_i-c_{ij}x_{i^{\prime}}x_{j^{\prime}}\in \Re, \ 1\leq i\neq j
\leq n; \;
\]
and

b) monomials:
\[
(\xi_i)^2 \in \Re^{\bot} , 1 \leq i \leq n.
\]
\end{definition}
\begin{remark} \cite{Maninpreprint}, (see also \cite{paul})
Note that if we set  $V=Span_{\textbf{k}}(x_1, x_2,  \cdots, x_n)$
$V^{*}= Span_{\textbf{k}}(\xi_1 ,\xi_2 \cdots, \xi_n),$ and define
a bilinear pairing $\langle \; \mid \; \rangle : V^{*}\otimes
V\longrightarrow \textbf{k}$ by $\langle \xi _i\mid x_j\rangle =
\delta_{ij},$ then the relations $\Re^{\bot}$ generate a subspace
in $V^{*}\otimes V^{*}$ which is orthogonal to the subspace of $
V\otimes V$ generated by $\Re.$
\end{remark}
\begin{definition}
\label{frobeniusalgebra}\cite{Maninpreprint}, \cite{Maninbook} A
graded algebra $A=\bigoplus_{i\ge0}A_i$ is called
\emph{a Frobenius algebra of dimension $d$},
(or \emph{a Frobenius quantum space of dimension $d$})
if
\begin{enumerate}
\item[(a)]
$dim(A_d)=1$,  $A_i =0,$ for $i>d;$
\item[(b)]
For all $j \geq 0$) the multiplicative map
$m: A_j\otimes A_{d -j} \rightarrow A_d$
is a perfect duality (nondegenerate pairing).

$A$ is called \emph{a quantum grassmann algebra} if in addition
\item[(c)]
$dim_ \textbf{k}A_i= \binom{d}{i} , \;\text{for}\; 1 \leq i \leq
d$ \end{enumerate}
\end{definition}
The following two theorems are  the main results of the paper.
\begin{theoremA}\label{frobeniustheoremF}
Let $X = \{x_1, \cdots, x_n\}$, let $\prec$ be the
degree-lexicographic order on $\langle X \rangle$. Suppose $\Fcal
=\textbf{k}\langle X \rangle/(\Re^{!})$ is a quadratic graded
algebra, which has precisely $\binom{n}{2}+n$ defining relations
\[
 \Re^{!} = \Re \bigcup \Re_1, \; \text{where}\;
\Re_1 = \{x_{j}x_{j}\}_{1\leq j\leq n},  \;
%\text{and  }\;
\Re=\{x_{j}x_{i} - c_{ij}x_{i'}x_{j'}\}_{1\leq i<j\leq n},
\]
and the set $\Re$  is such that:
\begin{enumerate}
\item[(a)] $\Re$ are relations of skew-polynomial type with
respect to $\prec$ (see Definition \ref{binomialringdef});
\item[(b)] $\Re$
is a Gr\"{o}bner basis for the ideal $(\Re)$ in $\textbf{k}\langle
X \rangle$.
\end{enumerate}
(In other words, $A = \textbf{k}\langle X \rangle/(\Re)$ is a
binomial skew polynomial ring).

Then
\begin{enumerate}
\item[(1)]
\label{frobeniustheoremF1} $\Re^{!}$ is a Gr\"{o}bner basis for the
ideal $(\Re^{!})$ in $\textbf{k}\langle X \rangle$ and the set of
monomials
\[
\Ncal^{!}= \{x_1^{\varepsilon_1}x_2^{\varepsilon
_2}\cdots x_n^{\varepsilon_n} \mid 0 \leq\varepsilon _i \leq 1, \
\text{for \ all} \ 1 \leq i \leq n \}
\]
is a $\textbf{k}$-basis of $\Fcal$.
\item[(2)]
$\Fcal$ is Koszul.
\item[(3)]
\label{frobeniustheoremF2} $\Fcal$  is a Frobenius algebra of
dimension $n$.
More precisely,  $\Fcal$ is graded (by length),
\begin{equation}\Fcal=\bigoplus_{i\ge0}\Fcal_i, \ \text{where}
\end{equation}
\begin{equation*}
\begin{split}
\Fcal_0 &= \textbf{k},\\
\Fcal_i& = Span_{\textbf{k}}\{u \mid u \in \Ncal^{!}
\;\text{and}\; \mid u\mid = i\}, \;
\text{for}\; 1\leq i \leq n,\\
\Fcal_n&= Span_{\textbf{k}}(W),  \; \text{where}\; W =  x_1x_2 \cdots x_n,\\
\Fcal_{n+j}&= 0 \; \text{for}\; j \geq 1.
\end{split}
\end{equation*}
\item[(
4)] \label{frobeniustheoremF3} Furthermore, $\Fcal$ is a
quantum grassmann algebra:
\begin{equation*}
\begin{split}
dim_ {\textbf{k}}\Fcal_i&= \binom{n}{i} , \;\text{for}\; 1 \leq i \leq n.\\
\end{split}
\end{equation*}
\end{enumerate}
\end{theoremA}

\begin{theoremB}
\label{maintheorem} Let $A= \textbf{k}\langle X \rangle/(\Re)$ be
a quantum binomial algebra. Then the following three conditions
are equivalent.
\begin{enumerate}
\item
\label{theoremB1} $A$ satisfies the weak cyclic condition. The
Koszul dual $A^{!}$ is Frobenius of dimension $n$, and has a
regular socle, see Definition \ref{regular1}.
%the
%monomial $W(\xi)$ generating the socle of $A^{!}$ has a regular
%presentation, $W= \xi_1 \cdots \xi_n,$
\item
\label{theoremB2} $A$ is a binomial skew polynomial ring, with
respect to some appropriate enumeration of $X.$
\item
\label{theoremB3} The automorphism $R=R(\Re):V^2 \longrightarrow
V^2$
%(see \ref{quantumquadraticmap}),
is a solution of the
classical Yang-Baxter equation, so $A$ is a Yang-Baxter algebra.
\end{enumerate}
Furthermore, each of these conditions implies that
\begin{enumerate}
\item[(a)]
\label{theoremBa}
There exists an enumeration of $X$, $X=\{x_1, \cdots, x_n \} $, such that the
set of ordered monomials $\Ncal _0$ forms a $\textbf{k}$-basis of $A$ , i.e.
$A$ satisfies an analogue of Poincar\'{e}-Birkhoff-Witt theorem;
\item[(b)] \label{theoremBb} $A$ is Koszul; \item[(c)]
\label{theoremBc} $A$ is left and right Noetherian.
\item[(d)]
\label{theoremBd}
$A$ is an Artin-Schelter regular domain.
\item[(e)]
\label{theoremBe} $A$ satisfies a polynomial identity.
\item[(f)]
\label{theoremBf} $A$ is catenary.
\end{enumerate}
\end{theoremB}

\section{The principal monomial and  regularity}

\begin{conventions}
\label{conventions1}
In this section we assume that
$A = A(\textbf{k}, X, \Re) = \textbf{k}\langle X\rangle/(\Re )$ is  a quantum binomial
algebra,
$\Scal _0 = \langle X ;\Re _0 \rangle$ is the associated quantum binomial semigroup.
$R: V^2 \longrightarrow V^2,$  and $r: X^2 \longrightarrow
X^2$, where  $R=R(\Re)$ and $r= r(\Re),$ are the maps associated
with $\Re,$ defined in \ref{quantumquadraticmap}. Furthermore,
till the end of the section we shall assume that the  Koszul dual
  $A^{!}$ is Frobenius.
  %see \ref{koszuldualdef}
\end{conventions}
\begin{remark}
\label{frobeniusremark} By our assumption

a) $A^{!}$ is graded by length:

\[
A^{!}=\bigoplus_{0\leq i \leq n}A_i^{!}, \ \text{where} \
dim(A_n^{!})=1;
\]
and

b) The multiplication function
$m: A_j^{!}\otimes A_{n -j}^{!} \rightarrow A_n^{!}$
is a non-degenerate pairing, for all $j \geq 0.$
\end{remark}
The one dimensional component $A_n^{!}$ is called \emph{the socle
of} $A^{!}$
\begin{notation}
For $m \geq 2,$ $\Delta_m = \{ x^m \mid x \in X\}$ denotes the
diagonal of $X^m.$
\end{notation}
\begin{definition}
\label{Oreconditions} Let $\Scal _0$ be a semigroup generated by
$X$. a) $\Scal _0$  satisfies \emph{the right Ore condition} if
for every pair $a, b\in X$ there exists a unique pair $x,y\in X$,
such that $ax=by$; b) $\Scal _0$  satisfies \emph{the left Ore
condition} if for every pair $a, b\in X$ there exists a unique
pair $z,t\in X$, such that $za=tb$.
\end{definition}
\begin{prooflemma} \ref{binomnchoose2relations}.
Suppose $A(\textbf{k}, X, \Re)$ is a quantum binomial algebra. By
Definition \ref{quantumquadraticdef} the relations in $\Re$ are
square-free, therefore $r(xx)=xx,$ and $\Lcal_x(x)=x= \Rcal_x(x)$
for every $x \in X.$ Suppose $x, y \in X, x \neq y$. The nondegeneracy condition implies
\[
\Lcal_x(y)\neq \Lcal_x(x)= x, \; \text{and} \;\Rcal_y(x)\neq
\Lcal_y(y)= y.
\]
It follows then that  the equality
\[
r(xy)= y^{\prime}x^{\prime} =\Lcal_x(y)\Rcal_y(x)
\]
implies
\begin{equation}
\label{eq1}
y^{\prime}\neq x,  x^{\prime}\neq y,
\end{equation}
therefore condition (c) holds. Clearly, (\ref{eq1}) implies  $r(xy)
\neq xy,$ so the relation $xy= y^{\prime}x^{\prime}$ belongs to
$\Re _0.$ It follows then that every monomial $xy \in X^2
\setminus \Delta_2$ occurs exactly once in $\Re _0,$ therefore in
$\Re,$ which verifies (a) and (b).  By \cite{T06}, Theorem 3.7, the
non-degeneracy of $r$, is equivalent to left and right Ore
conditions (see \ref{Oreconditions}) on the associated semigroup
$\Scal _0$ .
\end{prooflemma}

We recall some results which will be  used in the paper. The
following fact can be extracted from \cite{priddy}.
\begin{fact}
\label{priddykoszul} Suppose $A$ is a standard finitely presented
algebra with quadratic Gr\"{o}bner basis. Then $A$ is Koszul.
\end{fact}
\begin{theorem} (\cite{T06}, Theorem 9.7).
\label{skewimpliesYBE} Let  $A=\textbf{k} \langle X\rangle/(\Re)$
be a binomial skew-polynomial ring. Then  the automorphism
$R=R(\Re): V^2 \longrightarrow V^2,$ associated with $\Re,$ is a
solution of the Yang-Baxter equation, and $(X,r)$ is (a
square-free) set-theoretic solution of the Yang-Baxter equation.

Conversely, suppose $R: V^2 \longrightarrow V^2$ is a binomial
solution of the classical Yang-Baxter equation. Let $\Re= \Re(R)
\subset \textbf{k} \langle X\rangle$ be the quadratic binomial
relations defined via $R$. Then $X$ can be enumerated so, that the
Yang-Baxter algebra  $A=\textbf{k} \langle X\rangle/(\Re)$ is a
binomial skew polynomial ring. Furthermore every ordering $\prec$
on $X,$ $X = \{y_1, \cdots, y_n\},$ which makes the relations
$\Re$ to be of skew polynomial type, see Definition \ref{binomialringdef}
assures that $\Re$ is a Gr\"{o}bner basis with respect to $\prec$,
and the set of ordered monomials $\Ncal _{\prec}= \{y_1^{\alpha
_1}\cdots y_n^{\alpha _n}\mid \alpha _i \geq 0, 1 \leq i \leq n
\}$ is a $\textbf{k}$-basis for $A$.
\end{theorem}

%Let $\Xi = \{\xi_1, \cdots \xi_n \}$. We shall use the
%one-to-one correspondence
%\begin{equation}
%\iota :\langle \Xi \rangle \longrightarrow \langle X \rangle,
%\end{equation}
%where
%\begin{equation}
%\iota: \omega(\xi)=\xi_{i_1}\cdots \xi_{i_k} \mapsto x_{i_1}\cdots
%x_{i_k}= \omega .
%\end{equation}
For the following definition we do not assume Conventions \ref{conventions1} necessarily hold.
\begin{definition}
\label{translationdef} Let $\Xi = \{\xi_1, \cdots \xi_n \}$, be a
set of $n$ elements, which is  disjoint with $X.$  Let $T^{\xi}:
\langle X\rangle \longrightarrow \langle \Xi\rangle$, be the
semigroup isomorphism, extending the assignment $x_i \mapsto \xi
_i, 1 \leq i \leq n.$ If $\omega= \omega(x)= x_{i_1}\cdots x_{i_k}
\in \langle X\rangle,$ we call the monomial $T^{\xi}(\omega)=
\xi_{i_1}\cdots \xi_{i_k} \in \langle \Xi \rangle$
\emph{$\xi$-translation of} $\omega,$ and denote it by
$\omega(\xi).$ We define the
\emph{$\xi$-translation} of elements
$f \in \textbf{k} \langle X\rangle$, and of subsets
$\Re \subset\textbf{k} \langle X\rangle$ analogously,
and use notation $f(\xi)$ and  $\Re (\xi)$, respectively.
If $\Re _0 = \{\omega _i = \omega _i ^{\prime} \}_{i \in I}$ is a
set of semigroup relations in $\langle X\rangle $ by $\Re _0(\Xi)$
we denote the corresponding relations $\Re _0(\Xi): = \{\omega _i
(\xi)= \omega _i ^{\prime}(\xi) \}_{i \in I}$ in $\langle \Xi
\rangle$.
\end{definition}

Clearly the corresponding semigroups are isomorphic:
\[
\Scal _0=\langle X; \Re _0 \rangle \simeq \langle \Xi; \Re _0(\Xi)
\rangle
\]
and we shall often identify them. Let
\[
 \Scal ^{!}=\langle X;  \Re _0\bigcup \{(x_1)^2=0,
\cdots,  (x_n)^2=0\} \rangle
\]
Then the semigroup $\Scal ^{!}(\xi)$, associated with $A^{!}$, see Definition
\ref{Scalxi} is isomorphic to $\Scal ^{!}$.
\begin{definition}
\label{principaldef} Let $\Wcal= W(\xi) \in A^{!}$ be the monomial which
spans \emph{the socle}, $A_n^{!}$ of $A^{!}.$ Then the
corresponding monomial $W \in \Scal _0,$  is called \emph{the
principle monomial} of $A,$ we shall also refer to it as \emph{the
principle monomial} of $\Scal_0.$  A monomial $\omega \in \langle
X \rangle,$ is called \emph{a presentation of} $W$ if $W =
\omega$, as elements of $\Scal _0$.
\end{definition}
\begin{remark}
Clearly, $\mid W(\xi) \mid = n,$ so
$W(\xi)=\xi_{i_1}\xi_{i_2}\cdots \xi_{i_n},$ for some $i_j, 1\leq
i_j \leq n, 1 \leq j \leq n.$ Then the principal monomial
$W=x_{i_1}x_{i_2}\cdots x_{i_n} \in \langle X \rangle$, can be
considered as a monomial in $A,$ and in the semigroup $\Scal _0.$
Its
equivalence class (mod $\Re _0$) in $\langle X \rangle$
contains all monomials
$\omega \in
\langle X \rangle,$
which satisfy $\omega = W,$ in $\Scal _0.$ Clearly each such a
monomial $\omega$ has length $n,$ and is square-free. Furthermore,
$\omega = W,$ in $\Scal _0,$ if and only if $\omega (\xi) = c
W(\xi)$ in $A^{!},$ for an appropriate  $c\in \textbf{k}^{\times}.$
\end{remark}
We will define a special property of $W$, called \emph{regularity}
and will show that it is related to Artin-Schelter regularity of
$A.$ More precisely, for a quantum binomial algebra $A$ in which
the weak cyclic condition holds, the regularity of the principal
monomial $W$ implies Arin-Schelter regularity of $A$ and an
analogue of the Poincar\'{e}-Birkhoff-Witt theorem for $A$.

Till the end of the paper we shall often consider (at least) two
types  of equalities for monomials: a) $u=v$ as elements of $\Scal
_0 $ (or in $\Scal^{!}$) , and b) $u=v,$ as elements of the free
semigroup $\langle X \rangle.$ We remind that the equality a)
means that using the relations $\Re
_0$ (or the relations of
$\Scal^{!}$, respectively)  in  finitely many steps one can
transform $u$ into $v$ (and vice versa). The equality b) means
that $u$ and $v$ are equal as words (strings) in the alphabet $\{
x_1, \cdots , x_n\}.$ Clearly, b) implies a). To avoid ambiguity,
when necessary, we shall remind which kind of equality we
consider.
It follows from the Frobenius property
of $A^{!}$  that every $x_i$ , $1\leq i \leq n$ occurs as a head
(respectively, as a tail) of some presentation of $W.$

The presentation of the Koszul dual  $A^{!}$, in terms of
generators and relations is given in Definition \ref{koszuldualdef}.
\begin{definition}
\label{Scalxi} If we set $c_{xy}=1$ for all coefficients in the
defining relations of $A^{!}$, we obtain a new set of relations
which define a semigroup with zero. This way we associate
naturally to $A^{!}$, a semigroup with zero denoted by $\Scal(\xi)
^{!}$. As a set $\Scal(\xi) ^{!}$ is identified with the set
$\textbf{\Ncal}= Nor_{A^{!}}$ of normal monomials modulo the
(uniquely determined) reduced Gr\"{o}bner basis of $(\Re^{\bot}).$
Using the theory of Gr\"{o}bner basis it is easy to see that for
arbitrary $u, v \in \textbf{\Ncal}$ either

a) $uv = 0 $ in $A^{!}$, or

b) $uv = c w$ in $A^{!}$, with $c \in \textbf{k} ^{\times},$ and $w \in
\textbf{\Ncal},$ where the coefficient $c$ and the normal monomial
$w$ are uniquely determined, in addition $w\preceq uv$ in $\langle
X \rangle.$

We shall often identify $\Scal(\xi) ^{!}$ with the semigroup
$(\textbf{\Ncal}, \ast)$, where the operation $\ast$ is defined as
follows: for $u, v \in \textbf{\Ncal}$, $u\ast v: =0$ in case a)
and $u\ast v : = w$ in case b).
\end{definition}
\begin{remark} Note that $u*v = 0$ in $\Scal ^{!}(\xi),$ if and only if the
monomial $u(x)v(x),$ considered as a monomial in $\Scal _0$, has
some presentation, which contains a subword of the type $xx,$
where $x \in X.$ The shape of the defining relations of $A$, and
the assumption that $A^{!}$ is Frobenius, imply that a monomial $w
\in \langle X\rangle$ is a subword of some presentation of $W$,
($\omega = W$)  if and only if $w \neq 0$ as an element of $\Scal
^{!}.$
\end{remark}
\begin{definition}
\label{headsandtails} Let $w \in \Scal _0.$ We say that $h \in X$
is \emph{a head of $w$}  if $w$ can be presented (in $\Scal _0$)
as \[w = hw_1,\] where $w_1 \in \langle X \rangle$ is a monomial of
length
$\mid w_1 \mid = \mid w \mid -1$. Analogously, $t \in X$ is
\emph{a  tail of $w$} if \[w = w^{\prime}t\quad \text{(in $\Scal _0$)}\]  for
some $w^{\prime}\in \langle X \rangle,$ with $\mid w^{\prime} \mid
= \mid w \mid -1$.
\end{definition}
It follows from Remark \ref{frobeniusremark} b) that for every $i, 1\leq
i \leq n, $ there exists a monomial $\omega_i(\xi) \in \langle
\Xi\rangle , $ such that $\xi_i * \omega_i(\xi)= \Wcal$. Therefore
for every $i, 1 \leq i \leq n$, there exists a presentation
$W=x_i\omega _i,$ with $x_i$ as a head. Similarly, $x_i$ is a tail
of $W$ for every $i,1 \leq i \leq n. $ It is not difficult to
prove the following.
\begin{lemma}
\label{principalmonomial} The principal monomial $W$ of $\Scal _0$
satisfies the conditions:
\begin{enumerate}
\item $W$ is a monomial of length $n$.  There exist $n!$ distinct
words $\omega_i \in \langle X \rangle$,   $1 \leq i \leq n!$, for
which the equalities $\omega_i=W$  hold in $\Scal _0$. We call
them \emph{presentations} of $W$. \item Every $x\in X$ occurs as a
``head''  (respectively, as a ``tail'')  of some presentation of
$W$.
\[
W=x_1w_1^{\prime} =x_2w_2^{\prime}= \cdots x_nw_n^{\prime} \]
\[
W = \omega_1x_1=\omega_2x_2= \cdots \omega_nx_n.
\]
\item \label{W1} No presentation $\omega= W$, where  $\omega\in \langle X\rangle$
contains a subword of the form $xx,$ where $x \in X.$ \item
\label{W2} $W(\xi)$ spans the socle of the Koszul dual algebra
$A^!$. \item \label{W3} Every subword $a$ of length $k$ of
arbitrary presentation of $W$, has exactly $k$ distinct ``heads'',
$h_1, \cdots, h_k$, and exactly $k$ distinct ``tails'' $t_1,
\cdots, t_k.$ \item \label{W4} $W$ is the shortest monomial which
``encodes'' all the information about the relations $\Re _0,$
%that is the information about the solution  $(X, r)$.
More precisely, for any relation $(xy=y^{\prime}x^{\prime}) \in
\Re _0,$ there exists an $a\in \langle X\rangle$, such that
$W_1=xya$ and $W_2=y^{\prime}x^{\prime}a$  are (different)
presentations of $W$. \item \label{W5} If $W=ab$ is an equality in
$\Scal _0$, where $a,b \in \langle X \rangle,$  then there exists
a monomial $b^{\prime}\in \langle X\rangle,$ such that
$W=b^{\prime}a$ in $\Scal _0 $.
\end{enumerate}
\end{lemma}
%We shall see in Theorem \ref{regular2} that regularity of a special
%property of $W$ called regularity is essential for the properties
%of $A$, and $\Scal_0$.

Assume now that there exist a presentation
\begin{equation}
\label{regularpresentation} W= y_1y_2\cdots y_n,
\end{equation}
of $W$, in which  all  $y_1, y_2,  \cdots, y_n$ are pairwise
distinct, that is  $y_1, y_2, \cdots, y_n$ is a permutation of
$x_1, \cdots , x_n$. (The identity permutation is also allowed).
We  fix the degree-lexicographic order ``$\prec$'' on the free
semigroup $\langle y_1, \cdots, y_n \rangle = \langle X \rangle,$
assuming
\begin{equation}
\label{y} y_1 \prec y_2 \prec \cdots \prec y_n.
\end{equation}
We say that the  order $ \prec $ on $\langle  X \rangle $ is
\emph{associated with the presentation} \ref{regularpresentation}.

The theory of Gr\"obner bases, or the Diamond Lemma, see \cite{B}, implies that the set
of relations  $\Re _0$ determines a unique \emph{reduced Gr\"obner
basis} $\Gamma=\Gamma(\Re _0, \prec)$ in $\langle X \rangle$. In
general, $\Gamma$ is not necessarily finite. In fact, $\Re _0
\subseteq \Gamma,$ and every element of $\Gamma$ is of the form
$w=u,$ where  the monomials $u, w \in \langle X \rangle$ have
equal lengths $k \geq 2$, and $u \prec w.$ The monomial $w$ is
called the \emph{leading monomial} of the relation $w=u.$ (Note
that the relation $w=u$ follows from $\Re _0$, and holds in $\Scal
_0$.) A monomial $u \in \langle X \rangle$ is called \emph{normal
(mod $\Gamma$)}, if it does not contain as a subword any leading
monomial of some element of $\Gamma$. Clearly, if $u$ is normal,
then any subword $u^{\prime}$ of $u$ is normal as well. An
important property of the Gr\"obner basis $\Gamma$ is that every
monomial  $w \in \langle X \rangle$ can be reduced (by means of
reductions defined by $\Gamma$) to a  uniquely determined monomial
$ w_0 \in\langle X \rangle$, which is normal mod $\Gamma,$
and  such that $w=w_0$ is an equality in $\Scal_0  $. In addition
$w_0 \preceq w$ always holds in  $\langle X \rangle.$ The monomial
$w_0$ is called \emph{the normal form} of $w$ and denoted by
$Nor_{\Gamma}(w),$ or shortly $Nor(w)$.

Let $N=N(\Gamma)$ be the set of all normal (mod $\Gamma$)
monomials in $\langle X \rangle$. As a set $\Scal _0 $ can be
identified  with $N$. An  operation  ``$\ast$'' on  $N$ is
naturally defined as $u\ast v = Nor(uv)$, which makes $(N,\ast)$ a
semigroup, isomorphic to $\Scal _0 $.

It follows from the definition that there is an equality $\Re _0 =
\Gamma$ if and only if $\Scal _0$ is a semigroup of
skew-polynomial type (with respect to the ordering \ref{y}). The
Diamond lemma, \cite{B}, provides a recognizable necessary and
sufficient condition for   $\Re _0$ to be a Gr\"obner basis:
$\Re
_0$ is a Gr\"obner basis with respect to $\prec$, if and only if
every monomial of the shape $y_ky_jy_i$, with $n\geq k >j
> i \geq 1,$ can be reduced using $\Re _0$ to a uniquely
determined monomial of the shape $y_py_qy_r,$ with $p\leq q\leq
r.$
\begin{definition}
\label{regular1} Let
$W= W(r)$ be the
principal monomial of  $\Scal_0 $. We say that  $W= y_1y_2\cdots
y_n$, is a \emph{regular presentation} of $W$ if the following two
conditions are satisfied:

\begin{enumerate}
\item $y_1, y_2, \cdots, y_n$ is a permutation of $x_1, \cdots,
x_n$; and \item \label{d2} $y_1y_2\cdots y_n$ is the minimal
presentation of $W$ with respect to $\prec$ in $\langle X
\rangle$ (i.e. each $\omega \in \langle X\rangle$, such that
$\omega = W$ in $S$, satisfies $y_1y_2\cdots y_n\prec \omega$).

In this case we also say that $\prec$ is a \emph{regular order} in
$\langle X \rangle$

\end{enumerate}

We say that the Koszul dual $A^{!}$ has \emph{a regular socle}, if
the principal monomial $W$ has a regular presentation.
\end{definition}
\begin{remark}
\label{subwords} Let $W=y_1y_2\cdots y_n$ be a regular
presentation of $W$. It follows from the definition \ref{regular1}
that $Nor(W)=y_1y_2\cdots y_n$, or equivalently, the monomial
$y_1y_2\cdots y_n$ is normal (mod $\Gamma$.) Clearly, every subword of $y_1y_2 \cdots y_n$ is normal as well. In particular,  the
monomial $y_jy_{j+1}$ is normal for every $j, 1 \leq j <n$. Thus
$y_jy_{j+1}=zt \in \Re$ implies  $z\succ y_j, t \neq y_{j+1}.$
\end{remark}
\begin{example}
Let $X = \{x_1, x_2, x_3, x_4\}$, $S = \langle X; \Re_0 \rangle$ be
the
  semigroup with defining relations $\Re_0$:
\[\begin{array}{c}
x_1x_2=x_3x_4\quad x_1x_3=x_2x_4\quad x_4x_2=x_3x_1\\ x_4x_3=x_2x_1 \quad
x_1x_4=x_4x_1\quad x_2x_3=x_3x_2.
\end{array}
\]
Then the relations $\Re$ define a set-theoretic solution $(X, r)$
of the Yang-Baxter equation, therefore by \cite{T06}, $A^{!}$ is
Frobenius. Furthermore $x_1x_2x_3x_4=W$ is a presentation of $W$
as a product of pairwise distinct elements of $X$, but this
presentation is not regular. In fact, the monomial $x_3x_4$ is a
sub monomial of $W$, but it is not normal, since
$x_3x_4=x_1x_2$ (in $S$)  and $x_1x_2 \prec x_3x_4.$ Nevertheless  $W$ has
regular presentations. For example each of the monomials in the
following equalities gives a regular presentation of $W:$
$x_2x_3x_1x_4 = x_1x_4x_2x_3=x_4x_1x_3x_2=W$.
\end{example}
\begin{lemma}
\label{cancellation}  $\Scal^{!}$ has a cancelation  law on
nonzero products. More precisely, if $a, b, c  \in \Scal ^{!}$
then i) $ab= ac \neq 0$ implies $b=c$; ii) $ba = ca \neq 0$
implies $b=c.$
\end{lemma}
\begin{proof}
Conditions i) and ii) are analogous. We   shall prove i)  using
induction on the length $m$ of $a$. \begin{step} Let $\mid a\mid
=1,$ so $a \in X$.  Suppose for some monomials $b$ and $c$ one has:
\[
 ab = ac \neq 0.
\]
It follows then that $ab, ac$, and
therefore $b$ and $c$ are subwords of $W$. Clearly
$b$ and $c$
have equal lengths,
\[
\mid b\mid = \mid c \mid = k, \ k \geq 1.
\]
In the case when $k=1,$ the equality $ab=ac$ can not be a relation
because of the nondegeneracy property, therefore it is simply
equality of words in the free semigroup $\langle X \rangle,$ so
$b=c \in X$. Assume now that the length $k \geq 2,$ and,
\begin{equation}
\label{cancellation2} b \neq c.
\end{equation}
Note that  each of the monomials  $b$ and $c$ has
exactly $k$ heads, as a subword of $W$.  Let $H_b= \{b_1, \cdots, b_k \}$ be the set of
all heads of $b$ and $H_c = \{c_1, \cdots , c_k\}$ be the set of
heads of $c$. The inequality \ref{cancellation2} implies that
\begin{equation}
\label{cancellation3} H_b\neq H_c.
\end{equation}
The following relations hold in $\Scal_0,$ for appropriate
$b_i^{\prime}, c_i^{\prime}, a_i, a^{(i)} \in X,$ $1\leq i \leq
k.$
\begin{equation}\begin{array}{c}
\label{cancellation4} ab_i=b_i^{\prime}a_i, 1\leq i \leq
k, \\
 ac_i=c_i^{\prime}a^{(i)}, 1\leq i \leq k .
 \end{array}
\end{equation}
It follows from   (\ref{cancellation3}) and the non-degeneracy
property that there is an inequality of sets:
\[ \{b_i^{\prime} \mid 1 \leq i \leq k \} \neq
\{c_i^{\prime}\mid 1 \leq i \leq k\}.
\]
Clearly, then the set of heads of the
monomial $ab=ac$ is
\[
H_{ab} = \{ a\} \bigcup \{b_i^{\prime} \mid 1 \leq i \leq k \}\bigcup
\{c_i^{\prime}\mid  1 \leq i \leq k\}.
\]
By  the nondegeneracy condition one has $a \neq b_i^{\prime}$,
$a \neq c_i^{\prime},$ which together with (\ref{cancellation6})
imply that $H_{ab}$ contains at least $k+2$ elements.  This gives
a contradiction, since the monomial $ab$ is a subword  of $W$
therefore the number of its heads equals its length  $\mid ab \mid
= k+1$.
\end{step}
\begin{step} Assume the statement of the lemma is true for
all monomials $a, b, c$, with  $\mid a \mid \leq m.$ Suppose
$ab=ac\neq 0$ holds in $\Scal^{!}$, where $\mid a \mid = m+1.$ Let
$a= z_1 \cdots z_{m+1}$. Therefore $z_1*(z_2 \cdots z_{m+1}* b)=
z_1*(z_2 \cdots z_{m+1}* c),$ which by the inductive assumption
implies first that $(z_2 \cdots z_{m+1}* b)=(z_2 \cdots  z_{m+1}*
c)$, and again by the inductive assumption one has $b=c.$
\end{step}
\end{proof}
\begin{remark}
In some cases, when we study quadratic algebras, instead of
applying reductions to monomials of length $3$ ( in the sense of
Bergman \cite{B}), it is more convenient to study the action of
the infinite dihedral group, $\Dcal(\Re)$ generated by maps
associated with the quadratic relations,  as it is suggested
below.

Let $\Re$ be quantum binomial relations, $r = r(\Re)$ the
associated bijective map  $r: X^2 \longrightarrow X^2.$   Clearly
the two bijective maps $r_{ii+1}: X^3 \longrightarrow X^3$, $ 1
\leq i \leq 2$ , where $r_{12} = r\times Id_X$, and
$r_{23}=Id_X\times r$ are involutive. The infinite dihedral group,
\[
\Dcal= \Dcal(r)=\: _{\rm{gr}} \langle r_{12}, r_{23}: r_{12}^2 =
e, r_{23}^2=e \rangle
\]
acts on $X^3.$ The orbit $\Ocal_{\Dcal}(\omega)$ of $\omega \in
X^3$ consists of  all monomials $\omega^{\prime}\in X^3$ such that
$\omega^{\prime} = \omega$ is an equality in $\Scal _0.$ Clearly
each reduction $\rho$ applied to a monomial $\upsilon \in X^3$ can
be presented as $\rho(\upsilon)=r^{ii+1}(\upsilon),$ where $1\leq
i\leq 2.$ So every monomial $\omega^{\prime}$ obtained by a
sequence of reductions applied to $\omega$ belongs to
$\Ocal_{\Dcal}(\omega)$. The convenience of this approach is that
it does not depend on the enumeration of $X$ (therefore on the
chosen order $\prec$ on $\langle X\rangle.$
\end{remark}
\begin{lemma}
\label{weakcyclicconditionL} Suppose the quantum binomial algebra
$A= k\langle X; \Re\rangle$  satisfies the weak cyclic condition,
\ref{weakcycliccondition}. Let $\Ocal= \Ocal_{\Dcal}(\omega)$ be
an arbitrary  orbit of the action of $\Dcal$ on $X^3$. Then the
following conditions hold.
\begin{enumerate}
\item  $\Ocal \bigcap \Delta_3 \neq \emptyset$ if and only if
$\Ocal= \{xxx \},$ for some $x \in X$. \item $\Ocal \bigcap
((\Delta_2 \times X\bigcup X \times \Delta_2)\backslash \Delta_3))
\neq \emptyset$ if and only if $\mid \Ocal \mid = 3$.
\item
In each of the cases   $\omega =yyx,$ or $\omega =yxx,$
where $x, y \in X, x\neq y$,
the orbit $\Ocal_{\Dcal}(\omega)$ contains exactly $3$ elements.
More precisely,
if (by the weak cyclic condition) the following are relations in $\Scal _0:$
\[
yx= x_1y_1 ,  \;  \;  yx_1=x_2y_1 \;  \; \text{and} \; \;
y_1x_1=x_2y_2,
\]
%where $y_1=y_2$ or $x_1=x_2$ is also possible).
then there are equalities of sets:
\[\begin{array}{c} \Ocal_{\Dcal}(yyx)=\{ yyx, yx_1y_1,
x_2y_1y_1  \}\\
\Ocal_{\Dcal}(yxx)=\{ yxx, x_1y_1x, x_1x_1y_2  \} .
\end{array}\]
\end{enumerate}
Furthermore, suppose $\prec$ is an ordering on $X$ such that every
relation in $\Re _0$ is of the type $yx=x^{\prime}y^{\prime},$
where $y\succ x,$  $x^{\prime}\prec y^{\prime}$, and $y \succ
x^{\prime}.$ Then the orbit $\Ocal_{\Dcal}(y_1y_2y_3)$ with
$y_1\prec y_2 \prec y_3$ does not contain elements of the form
$xxy,$ or $xyy$, $x\neq y \in X$.
\end{lemma}

\begin{theorem}
\label{regular2} Let $A = A(\textbf{k}, X, \Re) = k\langle X; \Re
\rangle$ be a quantum binomial
algebra, let
$\Scal _0= \langle X ; \Re _0 \rangle$ be the associated semigroup,
and let $A^{!}$ be the Koszul dual of $A$.
We assume that the following conditions are satisfied:
  \begin{enumerate}
\item[(a)]
\label{regular22}
The weak cyclic condition is satisfied on $\Scal _0.$
\item[(b)]
\label{regular23}
The Koszul dual $A^{!}$ is Frobenius.
\item[(c)] \label{regular24} The principal monomial $W$ has a regular
presentation $W =y_1y_2\cdots y_n$ .
\end{enumerate}
Then
$\Scal _0=\langle y_1, y_2, \cdots, y_n; \Re _0\rangle$ is a
semigroup of skew-polynomial type (with respect to the order
$\prec$, where $y_1\prec y_2 \prec\cdots \prec y_n$ ). More precisely, the following conditions hold:
\begin{enumerate}
\item \label{t1} Each relation in $\Re _0,$ is of the form
$yz=z^{\prime}y^{\prime}$, where $y \succ  z$ implies
$z^{\prime}\prec y^{\prime}$, and $ y \succ z^{\prime}.$
%(and $z^{\prime}\succ y^{\prime}$ implies $y \prec z.$
\item \label{t2} The relations $\Re _0$ form a Gr\"obner
basis with respect to the  ordering $\prec$
on $\langle X \rangle.$
\item
\label{t3}
The relations $\Re$ form a Gr\"obner
basis with respect to the (degree-lexicographic) ordering $\prec$
on $\langle X \rangle,$ and $A$ is a binomial skew-polynomial ring.
  \item
  \label{t4}
The set of ordered monomials
\[
\Ncal= \{y_1^{\alpha_1} \cdots y_n^{\alpha_n}\mid \alpha_i \geq 0, 1 leq i \leq n \}
\]
forms a $\textbf{k}$-basis of $A$. That is $A$ is a PBW-type algebra.
\item
  \label{t5}
$A$ is Koszul. \item
  \label{t6}
$A$ is Artin-Schelter regular ring of global dimension $n$.
\end{enumerate}
\end{theorem}
We assume conditions  (a), (b) , (c) of the
hypothesis of the theorem are satisfied and prove two more
statements.
\begin{proposition}
\label{regular3} The following conditions hold in $\Scal _0$.
\begin{enumerate}
\item
\label{p1}
For any integer $j,$ $1 \leq j \leq n-1,$ there exists a unique
$\eta_j \in X,$ such that
\[
y_{j+1} \cdots y_n\eta_j = y_jy_{j+1}\cdots y_n.
\]
\item \label{p3} The elements $\eta_1, \eta_2, \cdots, \eta_{n-1}$
are pairwise distinct. \item \label{p4} For each $j$, $1 \leq j
\leq n-1,$ the set of heads $H_{W_j}$ of the monomial $W_j =
y_jy_{j+1}\cdots y_n$ is
\[
H_{W_j}= \{ y_j, y_{j+1},  \cdots , y_n \}.
\]
\item
\label{p5}
For any pair of integers $i,j$, $1 \leq i < j \leq n,$
the monomial $y_iy_j$ is normal. Furthermore, the unique
relation in which $y_iy_j$ occurs
has the form $y_{j^{\prime}}y_{i^{\prime}}= y_iy_j$, with
$j^{\prime} > i^{\prime},$ and $j^{\prime} > i.$
\end{enumerate}
\end{proposition}
\begin{lemma}
\label{regular4}
For each integer $j,$ $1 \leq j \leq n-1,$ let
$\xi_{j,j+1}, \cdots,  \xi_{j,n}$, $\eta_{j,j+1}, \cdots,  \eta_{j,n}$
be the elements of $X$ uniquelly determined by the relations
\begin{equation}
\begin{aligned}
\label{xij} &\xi_{j,j+1}\eta_{j,j+1}=y_jy_{j+1}\in \Re _0 \cr
&\xi_{j,j+2}\eta_{j,j+2}=\eta_{j,j+1}y_{j+2}\in \Re _0\cr &\cdots
\cdots \cdots\cr &\xi_{j,n-1}\eta_{j,n-1}=\eta_{j,n-2}y_{n-1} \in
\Re _0\cr &\xi_{j,n}\eta_{j,n}=\eta_{j, n-1}y_{n} \in \Re _0\cr
\end{aligned}
\end{equation}
Then for each $j$, $1 \leq j \leq n-1,$ the following conditions hold:
\begin{enumerate}
\item
\label{r1}
$\xi_{j,j+s}\neq \eta_{j,j+s-1},$ for all $s, 2\leq s \leq n-j$.
\item
\label{r2}
There is an equality of in $S$:
\[
\xi_{j,j+1}\xi_{j,j+2} \cdots \xi_{j,n} = y_{j+1} \cdots y_n.
%=y_{j+1}\vee y_{j+2}\vee \cdots \vee y_n.
\]
\item
\label{r3}
$y_{j+1}y_{j+2}\cdots y_n \eta_{j,n}= y_jy_{j+1}\cdots y_n$.
\item
\label{r4}
The elements $\eta_{j,n}, \eta_{j+1,n}, \cdots, \eta_{n-1,n}$ are
pairwise distinct.
\end{enumerate}
\end{lemma}
\begin{proof}
Condition (\ref{r1}) is obvious.
To prove the remaining conditions
we use decreasing induction on $j$, $1 \leq j \leq n-1.$

\textbf{Step 1.} $j=n-1.$ Clearly, $y_{n-1}y_n$ is normal thus
(see Remark \ref{subwords}) the relation in $\Re _0$ in which it occurs
has the shape $y_{n-1}y_n=\xi_{n-1,n}\eta_{n-1,n},$ with
$\xi_{n-1,n}\succ y_{n-1}.$ It follows then that $\xi_{n-1,n}=y_n$
and $y_{n-1}y_n=y_n\eta_{n-1,n},$, so $y_{n-1}y_n= y_{n-1}\vee
y_n.$. This gives (\ref{r2}). (\ref{r3}) and  (\ref{r4}) are clear.

\textbf{Step 2.} We first prove (\ref{r4}) for all $j, 1 \leq \j\leq
n-1.$ Assume that for all $k, \; n-1\geq k > j,$ the elements $y_k,
y_{k+1},\cdots, y_n$, $\xi_{k,k+1}, \cdots, \xi_{k,n}$,
$\eta_{k,k+1},\cdots, \eta_{k,n}$ satisfy
\begin{equation}
\begin{aligned}
\label{xik} &\xi_{k,k+1}\eta_{k,k+1}=y_ky_{k+1}\in \Re _0 \cr
&\xi_{k,k+2}\eta_{k,k+2}=\eta_{k,k+1}y_{k+2}\in \Re _0 \cr &\cdots
\cdots \cdots\cr &\xi_{k,n-1}\eta_{k,n-1}=\eta_{k,n-2}y_{n-1} \in
\Re _0 \cr &\xi_{k,n}\eta_{k,n}=\eta_{k, n-1}y_{n} \in \Re _0, \cr
\end{aligned}
\end{equation}
all $\eta_{j+1,n}, \eta_{j+2,n}, \cdots,\eta_{n-1,n}$
are pairwise distinct, and the modified
conditions
%(\ref{r2}),  (\ref{r3}), and
(\ref{r4}), in which ``$j$'' is replaced by ``$k$'' hold. Let
$\xi_{j,j+1}, \cdots \xi_{j,n}$, $\eta_{j,j+1}\cdots \eta_{j,n}$
satisfy (\ref{xik}). We shall prove  that $\eta_{j,n}\neq
\eta_{k,n},$ for all $k, j<k\neq n-1.$ Assume the contrary,
\[
\eta_{j,n}=\eta_{k,n}
\]
for some $k>j.$
Consider the relations
\begin{equation}
\label{e2} \xi_{j,n}\eta_{j,n}=\eta_{j,n-1}y_n, \quad
%\text{and} \
\xi_{k,n}\eta_{k,n}=\eta_{k,n-1}y_n.
\end{equation}
The Ore condition, (see Definition \ref{Oreconditions}),
%and  \ref{e1}
and (\ref{e2})
imply
\[
\label{e3} \eta_{j,n-1}=\eta_{k,n-1}.
\]
Using the same argument in $n-k$ steps we
obtain the equalities
\[
\eta_{j,n}=\eta_{k,n}, \eta_{j,n-1}=\eta_{k,n-1}, \cdots,  \eta_{j,k+1}=\eta_{k,k+1}.
\]
Now the relations
\[
\xi_{j,k+1}\eta_{j,k+1}=\eta_{j,k}y_{k+1},  \;
\xi_{k,k+1}\eta_{k,k+1}=y_ky_{k+1},
\]
and the Ore condition  imply $\eta_{j,k}= y_k.$
Thus, by (\ref{xik}) and (\ref{xij}) we obtain  a relation
\[\xi_{j,k}y_k=\xi_{j,k-1}y_k \in \Re_0.\] This is impossible, by
Lemma \ref{binomnchoose2relations} (iii). We have shown that the
assumption $\eta_{j,n}=\eta_{k,n}$, for some $k>j$, leads to a
contradiction. This proves (\ref{r4})
for all $j, 1 \leq j\leq n-1$.

We set
\begin{equation}
\label{eta}
\eta_1=\eta_{1,n}, \eta_2=\eta_{2,n}, \cdots, \eta_{n-1}=\eta_{n-1,n}.
\end{equation}

Next  we prove (\ref{r2}) and (\ref{r3}).

By the inductive assumption, for $k > j,$ we have
\[
\label{e4}
\xi_{k,k+1}\xi_{k,k+2}\cdots\xi_{k,n} = y_{k+1} \cdots y_n,
\]
and
\[
y_{k+1} \cdots y_n.\eta_{k+1}=y_k\cdots y_n.
\]
Applying the relations (\ref{xik}) one easily sees,  that
\[
\xi_{j,j+1}\xi_{j,j+2}\cdots\xi_{j,n}.\eta_{j,n}=y_jy_{j+1}\cdots y_n.
\]
Denote
\[
\omega_j=\xi_{j,j+1}\xi_{j,j+2}\cdots \xi_{j,n}.
\]
We have to show that the normal form,  $Nor(\omega_j)$, of
$\omega_j$ satisfies the equality of words
\[Nor(\omega_j) =
y_{j+1}y_{j+2} \cdots y_n, \quad\text{in}\; \langle X \rangle.\]
As a subword of length $n-j$ of the presentation $W=y_1y_2\cdots
y_{j-1}w_j\eta_{j,n}$, the monomial  $\omega_j$ has exactly $n-j$
heads
\begin{equation}
\label{e06}
h_1 \prec h_2 \prec \cdots \prec h_{n-j}.
\end{equation}
Since $Nor(\omega_j)=\omega_j$, is an equality  in $\Scal _0$, the
monomial $Nor(\omega_j)$ has the same heads as $\omega_j.$
Furthermore, there is an equality of words in $\langle X \rangle$,
  $Nor(\omega_j)=h_1\omega^{\prime}$, where $\omega^{\prime}$ is a
monomial of length $n- j -1.$ First we see that $h_1 \succeq y_{j}$
This follows immediately from the properties of the normal monomials
and the relations
\begin{equation}
\label{e7}
Nor(\omega_j)\eta_j=\omega_j\eta_j= y_jy_{j+1}\cdots y_n \in N .
\end{equation}
Next we claim that $h_1\succ y_j.$
Assume the contrary, $h_1=y_j.$
Then by (\ref{e7}) one has
\[
y_j\omega^{\prime}\eta_j=\omega_j\eta_j=y_jy_{j+1}\cdots y_n.
\]
The cancellation law  in $\Scal_0$ implies that
\[
\omega^{\prime}\eta_j=y_{j+1}\cdots y_n \in N.
\]
Thus $\eta_j$ is a tail of the monomial $y_{j+1}\cdots y_n.$ By
the inductive assumption, conditions (\ref{r2}) and (\ref{r3}) are
satisfied, which together with (\ref{eta}) give additional $n-j$
distinct tails of the monomial $y_{j+1}\cdots y_n$, namely
$\eta_{j+1}, \eta_{j+2}, \cdots \eta_{n-1}, y_n$. It follows then
that the monomial $y_{j+1}\cdots y_n$ of length $n-j$ has $n-j+1$
distinct tails, which is impossible. This implies  $h_1\succ y_j.$
Now since $\omega_j$ has precisely $n- j+1$ distinct heads, which
in addition satisfy (\ref{e06}) we obtain equality of sets
\[
\{h_1, h_2, \cdots , h_{n-j}\}=\{y_{j+1}, y_{j+2}, \cdots,  y_n\}.
\]
By the inductive assumption the heads of the monomial
$y_{j+1}y_{j+2}\cdots y_n$ are exactly $y_{j+1}, y_{j+2}, \cdots ,
y_n ,$ therefore,  there is an equality
\[\omega_j=
y_{j+1}y_{j+2}\cdots y_n  \quad \text{in  $\; \Scal _0$}.\] We have shown
(\ref{r3}). Now the equality
\[
y_{j+1}\cdots y_n\eta_j=y_{j}y_{j+1}\cdots y_n
\]
and the inductive assumption give that the heads of
$y_{j}y_{j+1}\cdots y_n$ are precisely the elements $y_j, y_{j+1},
\cdots, y_n.$ This proves (\ref{r2}). The lemma has been proved.
\end{proof}
\begin{proofproposition}
% \textbf{\ref{regular3}}.
Conditions  (\ref{p1}), (\ref{p3}), (\ref{p4}) of the proposition follow from
Lemma \ref{regular4}. We shall prove first that for any pair
$i,j,$ $1\leq i < j \leq n,$ the monomial $y_iy_j$ is normal.
Assume the contrary. Then there is a relation
\begin{equation}
\label{1p} (y_iy_j=y_{j^{\prime}}y_{i^{\prime}}) \in \Re_0,
\end{equation}
where
\[
y_{j^{\prime}}\prec y_i.
\]
Consider the monomial
\begin{equation}
\label{3p}
u=y_iy_j.y_{j+1}\cdots y_n\eta_{j-1} \eta_{j-2} \cdots \eta_{i+1}.
\end{equation}
We replace \ref{1p} in \ref{3p} and obtain
\[u=y_{j^{\prime}}y_{i^{\prime}}y_{j+1}\cdots
y_n\eta_{j-1}\eta_{j-2}\cdots \eta_{i+1},
\]
so $y_{j^{\prime}}$ is one of the heads of $u$. It follows from
\ref{regular4}.\ref{r3} that there is an equality in $\Scal_0$
$u=y_iy_{i+1}\cdots y_n = Nor(u)$. Since the inequality $Nor(u)
\preceq u$ always holds in $\langle X \rangle,$ $y_i$ is the
smallest ``head'' of $u$. But, by our assumption, the head
$y_{j^{\prime}}$ of $u$ satisfies $y_{j^{\prime}}\prec y_i$, which
gives a contradiction. We have proved that the monomial $y_iy_j$
is normal for every pair $i,j$, $1\leq i < j \leq n$.  Since the
number of relations is exactly $n\choose 2$ and each relation
contains exactly one normal monomial, this implies that all
monomials $x_jx_i,$ with  $1\leq i <j\leq n$, are not normal. It
follows then that each relation in $\Re _0$ has the shape
$y_jy_i=y_{i^{\prime}}y_{j^{\prime}},$ where $1 \leq i < j \leq
n,$ $1 \leq i^{\prime} < j^{\prime} \leq n,$ and $j > i^{\prime}$,
which proves (\ref{p4}) and (\ref{p5}).
\end{proofproposition}
\begin{lemma}
\label{rsolutionofYBE} The following conditions hold.
\begin{enumerate}
\item[(a)]
\label{rsolutionofYBEa}
The set of relations $\Re _0$  is Gr\"{o}bner basis with respect to the
ordering $\prec$ on $\langle X \rangle.$
\item[(b)]
\label{rsolutionofYBEb}
$\Scal_0$ is a semigroup of skew polynomial type.
\item[(c)]
\label{rsolutionofYBEc}
$(X, r)$ is a square-free solution of the set-theoretic
Yang-Baxter equation.
\item[(d)]
\label{rsolutionofYBEd}
$\Re$ is a Gr\"obner  basis of the ideal $(\Re)$
\item[(e)]
\label{rsolutionofYBEe}
$A$ is a binomial skew polynomial ring.
\item[(f)]
\label{rsolutionofYBEf}
The automorphism $R=R(\Re)$ is a solution of the classical
Yang-Baxter equation;
\end{enumerate}
\end{lemma}
\begin{proof}
We denote by $\Gamma$  the reduced Gr\"obner basis of the ideal $(\Re_0)$
and claim that $\Gamma =\Re_0.$
It will be enough to prove that the ambiguities $y_ky_jy_i,$ with  $k > j >i$,
do not give rise to new relations
in $\Scal _0.$ Or equivalently, the set $\Ncal _3$ of all monomials of length $3$
which are normal
modulo  $\Re _0$:
\[
\Ncal _3 = \{xyz\mid x,y,z\in X, \; \text{and} \; x\preceq y
\preceq z\},
\]
coincides with  $N_3= N\bigcap X^3$, the set of all monomials of
length $3$ which are  normal modulo $\Gamma $. Clearly, $N_3
\subseteq \Ncal _3.$

Let $\omega \in \Ncal _3.$ We have to show $Nor _{\Gamma}(\omega)
= \omega$. Four cases are possible:
\begin{equation}\begin{array}{ll}
\label{N1}
(i)\quad &\omega = y_iy_jy_k, 1 \leq i < j < k \leq n\\
 (ii)\quad &\omega = y_iy_iy_j, 1 \leq i < j  \leq n\\
(iii)\quad &\omega = y_iy_jy_j, 1 \leq i < j  \leq n\\
(iv) \quad &\omega = y_iy_iy_i, 1 \leq i  \leq n.
 \end{array}
\end{equation}
\textbf{Case 1.} Suppose \ref{N1} (i) holds. Assume, on the contrary, $\omega$
is not in $N_3.$ Then there is an  equality
\[\omega = y_iy_jy_k =
y_i^{\prime}y_j^{\prime}y_k^{\prime},\; \text{where} \; y_i^{\prime}\preceq
y_j^{\prime}\preceq y_k^{\prime},\] and  as elements of $\langle X
\rangle$, the two monomials satisfy
\begin{equation}
\label{N5}  y_i^{\prime}y_j^{\prime}y_k^{\prime}\prec y_iy_jy_k.
\end{equation}
By (\ref{N5}) one has
\begin{equation}
\label{N6}  y_i^{\prime}\preceq y_i .
\end{equation}
We claim that there is an inequality $y_i^{\prime}\prec
y_i.$ Indeed, it follows from Lemma \ref{weakcyclicconditionL}
that the orbit $\Ocal_{\Dcal}(y_iy_jy_k)$ does not contain
elements of the shape $xxy,$ or $xyy,$ therefore an assumption,
$y_i=y_i^{\prime}$ would imply $y_jy_k=y_j^{\prime}y_k^{\prime}$
with $y_j\prec y_k$ and $y_j^{\prime}\prec y_k^{\prime}$, which
contradicts Proposition \ref{regular3}. We have obtained that
$y_i^{\prime}\preceq y_i$.
One can easily see that there exists an $\omega \in \langle X
\rangle,$ such that
\[
(y_iy_jy_k) * \omega = y_i y_{i+1}\cdots y_n.
\]
The monomial $y_i y_{i+1}\cdots y_n$ is normal,  therefore the
normal form $(y_iy_jy_k) * \omega $ satisfies
\[
Nor(y_i^{\prime}y_j^{\prime}y_k^{\prime} * \omega) = Nor (y_iy_jy_k * \omega ) =y_i y_{i+1}\cdots y_n.
\]
Now the inequalities
\[
Nor(y_i^{\prime}y_j^{\prime}y_k^{\prime} * \omega) \preceq
y_i^{\prime}y_j^{\prime}y_k^{\prime} \prec y_i y_{i+1}\cdots y_n
\]
give a contradiction. It follows then that monomial $y_iy_jy_k,$
$i < j < k,$ is normal mod $\Gamma$.

\textbf{Case 2.} $\omega = y_iy_iy_j, 1 \leq i < j  \leq n$. It is not
difficult to see that the orbit $\Ocal = \Ocal_{\Dcal}(y_iy_iy_k)$
is the set
\[
\Ocal= \{ \omega= y_iy_iy_k, \omega_1=
y_iy_k^{\prime}y_i^{\prime},
\omega_3=y_k^{\prime\prime}y_i^{\prime}y_i^{\prime} \}
\]
where
\[
y_k^{\prime}y_i^{\prime}=y_iy_k \in \Re_0, \; \text{with} \; y_i\prec
y_k^{\prime}\succ y_i^{\prime}
\]
and
\[
y_k^{\prime\prime}y_i^{\prime}=y_iy_k^{\prime}\in \Re_0, \;
\text{with} \; y_k^{\prime\prime}\succ y_i^{\prime}.
\]
Therefore
\[
Nor( y_iy_iy_k) \in \Ocal_{\Dcal}(y_iy_iy_k)\bigcap \Ncal _3
=y_iy_iy_j
\]
We have shown that $Nor _{\Gamma} (\omega)= \omega.$

\textbf{Case 3} is analogous to Case 2. \textbf{Case 4.} is straightforward, since
all relations are square free. We have proved condition (a).

Condition (b) is straightforward.

We have shown that $\Scal_0= \langle X; \Re_0\rangle$ is a
semigroup of skew polynomial type. Clearly $r = r(\Re)= r(\Re
_0).$ It follows then from \cite{TM}, Theorem 1.1 that $(X,r)$ is
a solution of the set-theoretic Yang-Baxter equation which proves
(c).

We shall prove (d). It will be enough to show  that each ambiguity
$\omega =y_ky_jy_i,$ with $k>j>i$ is solvable.

Note first that since $(X,r)$ is a solution of the Yang-Baxter
equation, the group $\Dcal = $ is isomorphic to the dihedral group
$\Dcal _3$, and each monomial of length $3$ has an orbit
consisting either of $1,$ or $3$ or $6$ elements. Furthermore the
orbit  $\Ocal_{\Dcal}(\omega)$ consists of exactly $6$ elements.
This follows directly from lemma \ref{weakcyclicconditionL} (this
was proven first in \cite{TM}). Furthermore $\Ocal_{\Dcal}(\omega)$
contains exactly one ordered  monomial $\omega_0 =
y_{i_1}y_{j_1}y_{k_1},$ with $1 \leq i_1<j_1<k_1\leq n,$ which is
the normal form of $\omega$ (mod $\Re_0)$. Two cases are possible.
Either
\[ r^{12}r^{23}r^{12}(\omega) = \omega _0 =
r^{23}r^{12}r^{23}(\omega)
\]
or
\[ r^{12}r^{23}r^{12}r^{23}(\omega) = \omega _0 =
r^{23}r^{12}(\omega).
\]
Denote by $\Ocal_{\Re}(\omega)$ the set of all elements
$f\in A,$ which can be obtained by finite sequences of
reductions,
defined via the set of relations $\Re$
(in the sense of \cite{B}) applied to $\omega.$ In fact each
reduction $\rho$ applied to a monomial of length $3$, which is not
fixed under $\rho$ behaves as one of the automorphism $R^{12}$ and
$R^{23}$ but only in one direction, transforming each
monomial $\omega^{\prime}$ which is not ordered into
$\rho(\omega^{\prime}) = c_{pq}\omega^{\prime\prime}$, where
$c_{pq}$ is the coefficient occurring in the relation used for
$\rho$ and  $\omega^{\prime\prime}\prec \omega^{\prime}.$ So each
$f \in \Ocal_{\Re}(\omega)$ has the shape $f=c\omega^{\star},$
where $c \in \textbf{k}^{\times}$, and $\omega^{\star}$ is in the orbit
$\Ocal_{\Dcal}(\omega)$. We  know only that $ \Ocal_{\Re}(\omega)$
contains  $6$ elements, but the normal form, $\omega _0,$ might
occur with two distinct coefficients.

Assume now that the ambiguity $y_ky_jy_i$ is not solvable. Then
the orbit $\Ocal_{\Re}(\omega)$ contains two distinct elements
$c_1\omega _0$ and $c_2 \omega _0,$ with $c_1, c_2 \in
\textbf{k}^{\times},$ and $c_1, \neq c_2.$ On the other hand every $f \in
\Ocal_{\Re}(\omega)$ satisfies $f \equiv \omega$  (modulo $(\Re)$).
It follows then  $\omega _0 \in (\Re).$ One can find appropriate
$\eta_{s_1}, \cdots \eta _{s_{n-3}} \in X$ where $\eta _{j},$ $1
\leq j \leq n,$ are as in Proposition \ref{regular3} so that the
following equality holds in $\Scal_0$:
\[
y_{i_1}y_{j_1}y_{k_1} \eta_{s_1}\cdots \eta _{s_{n-3}}=W.
\]
But then there is an equality in $A$
\[
y_{i_1}y_{j_1}y_{k_1} \eta_{s_1}\cdots \eta _{s_{n-3}}= \alpha W
\]
for some  $\alpha \in \textbf{k}^{\times}$.  The element $\alpha W$
is in normal form therefore, $y_{i_1}y_{j_1}y_{k_1} =0$ in $A$
leads to a contradiction. It follows then that each ambiguity
$y_ky_jy_i,$ with $k>j>i,$ is solvable. Therefore  $\Re$ is a
Gr\"{o}bner basis, and  $A$ is a binomial skew polynomial ring.
This proves conditions (d) and (e). It follows from \cite{T06},
Theorem 9.7 (see also Theorem \ref{skewimpliesYBE}) that the automorphism
$R(\Re)$ is a solution of the classical Yang-Baxter equation.
\end{proof}
\begin{proofoftheorem}
%\textbf{\ref{regular2}}.
Condition (\ref{t1}) follows from Proposition \ref{regular3}.
(\ref{p5}). Lemma \ref{rsolutionofYBE}  implies (\ref{t2}) and
(\ref{t3}). Clearly (\ref{t3}) implies (\ref{t4}). It is known
that every standard finitely presented algebra with quadratic
Gr\"{o}bner basis is Koszul, Fact  \ref{priddykoszul}, which
implies (5). We have already shown that $A$ is a binomial skew
polynomial ring. It follows then from the proof of Theorem
\ref{Frobeniusth} that $A$ has global dimension $n$. Now a
result of  Stafford and Zhang, \cite{TobyZ}
see also P. Smith's,
see \ref{factkoszul},  implies that $A$ is
Gorenstein and therefore, $A$ is Artin-Schelter regular.
\end{proofoftheorem}

\section{The koszul dual of a binomial skew polynomial ring is  Frobenius }

In this section we  study the Koszul dual $A^{!}$ of a binomial
skew polynomial ring $A.$ We prove Theorem A,
%\ref{frobeniustheoremF}
which guarantees  Frobenius property for a
class of  quadratic algebras with specific relations. This class
includes the Koszul dual $A^{!}$. The main result of the section,
\ref{Frobeniusth}, shows that the binomial skew polynomial rings
with $n$ generators provide a class of Artin-Schelter regular
rings of global dimension $n$. The first proof of this theorem
(1995)  was given in \cite{T3},  where we used combinatorial
methods to show that $A^{!}$ is Frobenius, and then a result of
P. Smith \label{factkoszul}, to show that $A$ is regular. In
\cite{TM} this result was improved by a different argument, which
uses the good algebraic and homological properties of semigroups
of $I$-type to show that $A$ is an Artin-Schelter regular domain.
We  present here the original combinatorial proof of the Frobenius
properties of $A^{!}$, which has not been published yet and uses a
technique which might be useful in other cases of (standard)
finitely presented algebras.
\begin{theorem}
\label{Frobeniusth} Let $A$ be a skew-polynomial ring with
binomial relations. Then
\begin{enumerate}
\item
\label{Frobeniusth1} The Koszul dual $A^{!}$ is Frobenius.
\item
\label{Frobeniusth2} $A$ is Artin-Schelter regular ring of global
dimension $n$.
\end{enumerate}
\end{theorem}

Our proof is combinatorial, we deduce the Frobenius property of $A^{!}$
from its defining relations. We
use Gr\"obner  basis techniques, the cyclic condition in $A$, and
study more precisely  the computations in the associated semigroup
$\Scal^{!}$. Next we recall the following result which will be used
to deduce the Gorenstein property of $A$.
\begin{proposition}
\cite{paul}, Proposition 5.10.
\label{factkoszul}
Let $A$ be a Koszul
algebra of finite global dimension. Then $A$ is Gorenstein if and
only if $A^{!}$ is Frobenius. \end{proposition}

We keep the notation from the previous sections. As before we
denote the set of generators of $A^{!}$ as $\Xi= \{\xi_1, \xi_2,
\cdots , \xi_n.\}$
\begin{remark}
\label{skewimpliescycliccondition} In \cite{T96}, Theorem 3.16 (see
also \cite{T1}) was shown that every binomial skew polynomial ring
$A$ satisfies \emph{the cyclic condition}, a condition stronger
than the weak cyclic condition , see Definition \ref{weakcycliccondition}. So the
algebra $A$, satisfies the conditions of Definition \ref{weakcycliccondition}.
One can easily
deduce from the relations of $A^{!}$, see Notation \ref{allrelations} that
it also satisfies the conditions of Definition \ref{weakcycliccondition}.
\end{remark}
We need the
explicit relations of $A^{!}$.

Let  $A=\textbf{k} \langle X\rangle/(\Re)$ be a binomial
skew-polynomial ring, with a set of relations
\begin{equation}
\label{relationsre} \Re=\{ x_{j}x_{i} -
c_{ij}x_{i^\prime}x_{j^\prime}\}_{1\leq i<j\leq n}.
\end{equation}
where for each pair $1 \leq i < j \leq n,$ the relation
$x_{j}x_{i} - c_{ij}x_{i^\prime}x_{j^\prime},$ satisfies $j>
i^{\prime}$, $i'< j',$ and  $c_{ij} \in \textbf{k}^{\times}.$
Furthermore, the relations $\Re$ form a Gr\"obner basis, with
respect to the degree-lexicographic order on $\langle X\rangle.$

\begin{notation}
\label{allrelations} Let $\Xi=\{ \xi_1, \cdots , \xi_n\}$ be a set
of indeterminates, $\Xi\bigcap X = \emptyset.$ Consider the
following subsets of the free associative algebra $\textbf{k}
\langle \Xi\rangle$:
\[ \Re^{*}=\{ \xi_{j}\xi_{i} +
(c_{ij})^{-1}\xi_{i^\prime}\xi_{j^\prime}\}_{1\leq i<j\leq n}.
\]
We call $\Re^{*}$ \emph{the dual relations, associated to} $ \Re. $ Let
\[\begin{array}{llclll}
\Re _1&=&\{  (x_j)^2\}_{1\leq j\leq n},\quad&  \Re^{!}&=& \Re \bigcup \Re_1, \\
&&&&&\\
\Re _1^{*}&=&\{ (\xi _j)^2 \}_{1\leq j\leq n},
\quad&
 \Re^{\bot}&=& \Re^{*}\bigcup \Re_1^{*}.
\end{array}
\]
\end{notation}
It follows from the definition of Koszul dual \ref{koszuldualdef}
that:
\begin{remark}
Let $A=\textbf{k} \langle X\rangle/(\Re)$ be a binomial
skew-polynomial ring, with a set of relations $\Re$ as in
(\ref{relationsre}). Then the Koszul dual $A^{!}$ has the following
presentation via generators and relations:
\begin{equation}
\label{A!} A^{!}=\textbf{k} \langle \Xi \rangle /(\Re^{\bot}).
\end{equation}
\end{remark}

%\begin{definition}
%\label{translationdef} Let  $T^{\xi}: \langle X\rangle
%\longrightarrow \langle \Xi\rangle$, be the semigroup isomorphism,
%extending the assignment $x_i \mapsto \xi _i, 1 \leq i \leq n.$ If
%$\omega= \omega(x)= x_{i_1}\cdots x_{i_k} \in \langle X\rangle,$
%we call the monomial $T^{\xi}(\omega)= \xi_{i_1}\cdots \xi_{i_k}
%\in \langle \Xi \rangle$  \emph{$\xi$-translation of} $\omega,$
%and denote it by
% $\omega(\xi).$ We define the
%\emph{$\xi$-translation} of elements
% $f \in \textbf{k} \langle X\rangle$, and of subsets
% $\Re \subset\textbf{k} \langle X\rangle$ analogously,
% and use notation $f(\xi)$ and  $\Re (\xi)$, respectively.
%\end{definition}

The next lemma is straightforward
\begin{lemma}
\label{translationlemma} Let $\omega\in \langle X\rangle.$ Suppose
$\Re \subset \textbf{k} \langle X\rangle$ is a set of quantum
binomial relations, and $\Re^{*} \subset \textbf{k} \langle
\Xi\rangle,$ is the associated dual relation set. Let $\Re _0$ and
$\Re^{*}_0$ , respectively, be the semigroup relations associated
with $\Re$ and $\Re^{*}, see Definition \ref{associatedsemigroup}.$ Then the
following conditions hold:
\begin{enumerate}
\item There is an equality $(\omega(\xi))(x)= \omega.$
\item
$
(\Re _0)(\xi)= (\Re^{*})_0= (\Re_0)^{*}.
$
\item The $\xi$-translation isomorphism $T^{\xi}$ induces
(semigroup) isomorphisms:

a) between the associated semigroups
\[
\Scal_0=\Scal_0(X, \Re_0)= \langle X ; \Re_0\rangle \simeq
\Scal_0(\xi)=\Scal _0(\Xi, \Re^{*}_0)= \langle \Xi ; \Re^{*}_0 ,
\rangle
\]
and

b) between the "Koszul-type" semigroups
\[
\Scal^{!}= \langle X ; \Re_0 \bigcup \Re_1\rangle \simeq
(\Scal(\xi))^{!}=\langle \Xi ; \Re^{*}_0  \bigcup \Re_1^{*}
\rangle.
\]
\end{enumerate}
\end{lemma}
For our purposes it will be often more convenient to perform
computations and arguments in  $\Scal_0,$  $\Scal^{!}$ and $A$,
respectively, and then "translate" the results for $\Scal_0(\xi),$
$(\Scal(\xi))^{!}$ and $A^{!}.$
\begin{lemma}
\label{frobeniuslemma1}
In notation \ref{allrelations},
the following conditions are equivallent:
\label{gbases}
\begin{enumerate}
\item
\label{reX} $\Re$ is a Gr\"obner basis of the ideal $(\Re)$ in
$\textbf{k} \langle X\rangle.$
\item
\label{re*} $\Re^{*}$ is a Gr\"obner basis of the ideal
$(\Re^{*})$ in $\textbf{k} \langle \Xi \rangle.$
\item
\label{re!} $\Re^{!}$ is a Gr\"obner basis of the ideal
$(\Re^{!})$ in $\textbf{k} \langle X \rangle.$
\item
\label{rebot} $\Re^{\bot}$ is a Gr\"obner basis of the ideal
$(\Re^{\bot})$ in $\textbf{k} \langle \Xi \rangle.$
\end{enumerate}
\end{lemma}
\begin{proof}
Let $V = Span X,  \;  V^{*}=Span \ \Xi.$

We show first the implication $\ref{reX} \Longrightarrow
\ref{re*}.$ The implication $\ref{re*} \Longrightarrow \ref{reX}$
is analogous.

Suppose condition \ref{reX} holds. Clearly, this
implies that the algebra $A(\textbf{k}, X, \Re)$ is a binomial
skew polynomial ring. It follows then from Theorem \ref{skewimpliesYBE}
that the automorphism $R= R(\Re): V^2 \longrightarrow V^2$ is a
solution of the Yang-Baxter equation. It is not difficult to see
that $R^{*}= R(\Re^{*}): (V^{*})^2 \longrightarrow (V^{*})^2$ is
also a solution of the Yang-Baxter equation. Clearly the relations
$\Re^{*}$  are of skew-polynomial type. Hence by
theorem \ref{skewimpliesYBE}, $\Re^{*}$ is a Gr\"obner basis
of the ideal $(\Re^{*})$ in $\textbf{k} \langle \Xi \rangle$.

The implication $\ref{reX} \Longrightarrow \ref{re!}$  is verified
directly by Gr\"obner bases technique,
that is one shows that all
ambiguities are solvable, see the Diamond Lemma, \cite{B}. Clearly
there are three types of ambiguities: a) $x_kx_jx_i, n \geq k
> j
> i \geq 1$ , b) $x_jx_ix_i, n \geq j>i \geq 1,$ and
c) $x_jx_jx_i, n \geq j>i \geq 1.$ All ambiguities of the type a)
are solvable, since by (\ref{reX})  $\Re$ is a Gr\"{o}bner basis.
We will show that all ambiguities of type b) are solvable. Let
$j,i$ be a pair of integers, with $n \geq j > i \geq 1$. Consider
the ambiguity $x_jx_ix_i$. It follows from the cyclic condition
\ref{weakcycliccondition} that there exist integers $i_1, j_1,
j_2$, with $1 \leq i_1 < j_1, j_2 \leq n$ such that $\Re$ contains
the relations: $x_jx_i - c_{ij}x_{i_1}x_{j_1}$ and $x_{j_1}x_i
-c_{ij_1}x_{i_1}x_{j_2}$, where $c_{ij}$ and $c_{ij_1}$ are
nonzero coefficients. This gives the following sequence of
reductions:
\begin{equation*}
x_jx_ix_i \longrightarrow ^{R^{12}}(c_{ij}x_{i_1}x_{j_1})x_i
\longrightarrow ^{R^{23}}
c_{ij}x_{i_1}(c_{ij_1}x_{i_1}x_{j_2}) =
c_{ij}c_{ij_1}[x_{i_1}x_{i_1}]x_{j_2} \longrightarrow ^{R^{12}} 0.
\end{equation*}
The other possible way of reducing $x_jx_ix_i$ is
\begin{equation*}
x_jx_ix_i \longrightarrow^{R^{23}} 0.
\end{equation*}
We have proved that all ambiguities of the type b) are solvable.
An analogous argument shows that the ambiguities of the type c)
are also solvable. Thus $\Re^{!}$ is a Gr\"{o}bner  basis of the
ideal $(\Re^{!})$ in $\textbf{k} \langle X \rangle.$
\end{proof}

\begin{corollary}
\label{frobeniuscorrolary2} Let $A=\textbf{k} \langle
X\rangle/(\Re)$ be a binomial skew-polynomial ring, let  $A^{!}$
be its Koszul dual. Let $\Fcal=\textbf{k} \langle X\rangle/(\Re
^{!})$ Then
\begin{enumerate}
\item
$\Fcal$ has a $\textbf{k}$-basis the set
\begin{equation*}
\label{Ncal!} \Ncal^{!}= \{x_1^{\varepsilon_1}x_2^{\varepsilon
_2}\cdots x_n^{\varepsilon_n} \mid \varepsilon _i = 0, 1, \
\text{for \ all} \ 1 \leq i \leq n \}.
\end{equation*}
\item
$A^{!}$ has a $\textbf{k}$-basis the set
\begin{equation*}
\label{Ncalxi!} \Ncal(\xi)^{!}=
\{\xi_1^{\varepsilon_1}\xi_2^{\varepsilon _2}\cdots
\xi_n^{\varepsilon_n} \mid \varepsilon _i = 0, 1, \ \text{for \
all} \ 1 \leq i \leq n \}.
\end{equation*}
\item
The principal monomial of $A$ has a regular presentation $W=x_1
x_2 \cdots x_n.$
\item
The socle of $A^{!}$  is one dimensional and is generated by the monomial
$W(\xi) =\xi_1 \xi_2 \cdots \xi_n.$
\end{enumerate}
\end{corollary}
\begin{remark}
\label{Scal!} The semigroup $\Scal^{!}= \langle X ; (\Re_0 \bigcup
\Re_1)\rangle$ can be presented as $\Scal^{!}\simeq \Scal _0 /
(\Re_1).$ It is a semigroup with $0$, $xx= 0$ for every $x \in X.$
To make the computations in $\Scal^{!}$ we compute modulo the
relations $\Re _0,$ and keep in mind that $\omega \in \langle X
\rangle,$ is equal to $0$ in $\Scal^{!}$  if and only if it can be
presented as $\omega= \omega^{\prime}$ in $\Scal _0 ,$ where
$\omega^{\prime}= axxb \in \langle X \rangle,$ for some $x \in X,
a, b \in \langle X \rangle.$ Denote
\begin{equation*}
\Ncal^{!}_0= \Ncal^{!}\bigcup \{0\}
\end{equation*}
We can identify $\Scal^{!}$
with the semigroup $(\Ncal^{!}_0, \ast)$ where the operation
$\ast$ on $\Ncal^{!}_0$ is defined as follows: for $u, v \in
\Ncal^{!}_0$, either a) $u\ast v=0$ and this is true if and only
if the normal form $Nor_{\Re_0}(uv)$ contains some subword of the
shape $xx, x \in X$, or  b)  $u\ast v = w \in \Ncal^{!}$, where
$Nor_{\Re_0}(uv)=w$ (or equivalently $Nor_{\Re}(uv)= c w, $ for
some nonzero coefficient $c$).

All relations in $\Scal _0,$ which do not involve subwords of the
shape $xx$ are preserved in $\Scal^{!}$. In particular,  the
cyclic conditions are in force.

If $u, v, w \in \Ncal _0$ and $u \ast  w \neq 0$ (that is $u \ast
w \in \Ncal^{!}$), then each of the equalities $u \ast w = v \ast
w$ and $w \ast u = w \ast  v$ implies $u = v$, i.e. $(\Ncal^{!}_0,
\ast)$ has cancellation low for non-zero products.
\end{remark}

Theorem A
%\ref{frobeniustheoremF}
verifies the Frobenius property
for each quadratic algebra with relations of the type $\Re^{!}.$
We prove first some more statements under the hypothesis of
Theorem A.

Before proving the theorem we need some more statements.

Clearly the assumption that $A=\textbf{k} \langle X\rangle/(\Re)$
is a binomial skew-polynomial ring,  implies that $\Scal _0
=\langle X;\Re _0 \rangle$ is a semigroup of skew-polynomial type.
(with respect to the degree-lexicographic order $<$ on $\langle X
\rangle$ defined by $x_1 < x_2 < \cdots < x_n$. It is proven in
\cite{T96}, that  $\Scal _0$ satisfies the cyclic condition,
therefore Ore condition holds.
Furthermore $\Scal _0$ is with cancellation law, \cite{T3}.
Proposition \ref{regular5} is true for arbitrary semigroup of
skew-polynomial type. In some parts we use argument similar to the
proof of Proposition \ref{regular3}, but  we prefer to give sketch of the
proofs explicitly, since they  are made under different
assumptions.

\begin{proposition}
\label{regular5} Let $\Scal _0 =\langle X;\Re _0 \rangle$ be a
semigroup of skew-polynomial type, with respect to the
degree-lexicographic ordering $\prec$ on $\langle X \rangle$. Then
the following conditions are satisfied.
\begin{enumerate}
\item \label{regular51} The monomial $W_1 =x_1x_2\cdots x_n$ is
normal. \item \label{regular52} For any $j, 1\leq j \leq n-1$,
there exist a unique $\eta_j\in X$, such that $x_{j+1}\cdots
x_n\eta_j=x_jx_{j+1}\cdots x_n.$ \item \label{regular53} The
elements $\eta_1, \cdots, \eta_{n-1}$ are pairwise distinct. \item
\label{regular54} For every $j, 1\leq j \leq n-1,$ the monomial
$W_j=x_jx_{j+1}\cdots x_n$ has exactly $n-j+1$ heads, namely :
\begin{equation*}
H_{W_j}=\{x_j, x_{j+1},\cdots, x_n \}.
\end{equation*}
\item
\label{regular55} For any $j, 1\leq j \leq n-1$, there exist a
unique $\theta_{j+1}\in X$, such that $\theta_{j+1}x_{1}\cdots
x_j=x_1x_{2}\cdots x_{j+1}.$
\item
\label{regular56} The elements $\theta_2, \cdots, \theta_{n}$ are
pairwise distinct.
\item
\label{regular57} For every $j, 1\leq j \leq n-1$ the monomial
$\omega_j=x_1x_2 \cdots x_j$ has exactly $j$ tails, namely
%its set of tails $T_{\omega_j}$ satisfies:
\begin{equation*}
T_{\omega _j}=\{x_1, x_2,\cdots, x_j \}.
\end{equation*}
In particular, every $x_i, 1 \leq i \leq n$ occurs as a head and
as a  tail of the monomial $W_1=x_1x_2\cdots x_n = \omega _n$,
\item
\label{regular58} The monomial  $W_1$ is the principal monomial of
$\Scal _0$ with a regular presentation $W_1=x_1x_2 \cdots x_n.$
\end{enumerate}
\end{proposition}
Under the assumption of Proposition \ref{regular5} we prove first
Lemma \ref{regular6}.  Although the statements of Lemmas \ref{regular6} and
\ref{regular4}, are analogous, due to the different hypotheses, we need
different arguments for their proofs.
\begin{lemma}
\label{regular6} For each integer $j,$ $1 \leq j \leq n-1,$ let
$\zeta_{j,j+1}, \cdots,  \zeta_{j,n}$, $\eta_{j,j+1}, \cdots,
\eta_{j,n}$ be the elements of $X$ uniquely determined by the
relations
\begin{equation}
\begin{aligned}
\label{zetaj} &(\zeta_{j,j+1}\eta_{j,j+1}=x_jx_{j+1})\in \Re _0 \cr
&(\zeta_{j,j+2}\eta_{j,j+2}=\eta_{j,j+1}x_{j+2})\in \Re _0\cr
&\cdots \cdots \cdots\cr
&(\zeta_{j,n-1}\eta_{j,n-1}=\eta_{j,n-2}x_{n-1}) \in \Re _0\cr
&(\zeta_{j,n}\eta_{j,n}=\eta_{j, n-1}x_{n}) \in \Re _0.\cr
\end{aligned}
\end{equation}
Then for each $j$, $1 \leq j \leq n-1,$ the following conditions
hold:
\begin{enumerate}
\item \label{r61} $\zeta_{j,j+s}\neq \eta_{j,j+s-1},$ for all $s,
2\leq s \leq n-j$.

\item \label{r62} The following  are equalities
in $\Scal _0$:
\[\begin{array}{lcl}
\zeta_{j,j+1}\zeta_{j,j+2} \cdots \zeta_{j,n}&=& x_{j+1} \cdots x_n
%=y_{j+1}\vee y_{j+2}\vee \cdots \vee y_n\\
\\&&\\
x_{j+1}x_{j+2}\cdots x_n \eta_{j,n}&=& x_jx_{j+1}\cdots
x_n.
\\&&\\
\end{array}
\]
%\item
%\label{r63} $x_{j+1}x_{j+2}\cdots x_n \eta_{j,n}= x_jx_{j+1}\cdots
%x_n$.
\item
\label{r64} The elements $\eta_{j,n}, \eta_{j+1,n}, \cdots,
\eta_{n-1,n}$ are pairwise distinct.
\end{enumerate}
\end{lemma}
\begin{proof}
Condition \ref{r61} is obvious. To prove the remaining conditions
we use decreasing induction on $j$, $1 \leq j \leq n-1.$

\textbf{Step 1.} $j=n-1.$ Clearly, $x_{n-1}x_n$ is normal thus
(cf. Remark \ref{subwords}) the relation in $\Re_0$ in which it occurs has
the shape $x_{n-1}x_n=\zeta_{n-1,n}\eta_{n-1,n},$ with
$\zeta_{n-1,n}\succ x_{n-1}.$ It is clear  then that
$\zeta_{n-1,n}=x_n$ and $x_{n-1}x_n=x_n\eta _{n-1,n}.$ Hence the
set of heads of $x_{n-1}x_n$ is $ \{ x_{n-1},  x_n \}$. This gives
\ref{r62},    \ref{r64} is trivial.

\textbf{Step 2.} Using decreasing induction on $j$ we first prove condition (\ref{r64}) for all $j, 1 \leq
\j\leq n-1$. (Step 1, $j=n-1$ gives the base for the induction.
 Assume that for all $k, n-1\geq k > j,$ the elements
$x_k, x_{k+1},\cdots, x_n$, $\zeta_{k,k+1}, \cdots, \zeta_{k,n}$,
$\eta_{k,k+1},\cdots, \eta_{k,n}$ satisfy
\begin{equation}
\begin{aligned}
\label{zetak} &(\zeta_{k,k+1}\eta_{k,k+1}=x_kx_{k+1})\in \Re _0 \cr
&(\zeta_{k,k+2}\eta_{k,k+2}=\eta_{k,k+1}x_{k+2})\in \Re _0\cr
&\cdots \cdots \cdots\cr
&(\zeta_{k,n-1}\eta_{k,n-1}=\eta_{k,n-2}x_{n-1}) \in \Re _0\cr
&(\zeta_{k,n}\eta_{k,n}=\eta_{k, n-1}x_{n}) \in \Re _0; \cr
\end{aligned}
\end{equation}
all $\eta_{j+1,n}, \eta_{j+2,n}, \cdots,\eta_{n-1,n}$ are pairwise
distinct, and the modified conditions
%(\ref{r2}),  (\ref{r3}), and
(\ref{r64}), in which ``$j$'' is replaced by ``$k$'' hold. Let
$\zeta_{j,j+1}, \cdots \zeta_{j,n}$, $\eta_{j,j+1}\cdots
\eta_{j,n}$ satisfy (\ref{zetak}). We shall prove  that
$\eta_{j,n}\neq \eta_{k,n},$ for all $k, j<k\neq n-1.$ Assume the
contrary,
\begin{equation}
\label{e61} \eta_{j,n}=\eta_{k,n}
\end{equation}
for some $k>j.$ It follows from (\ref{e61}), the relations
\begin{equation*}
\label{e62} \xi_{j,n}\eta_{j,n}=\eta_{j,n-1}y_n, \; \text{and} \;
\xi_{k,n}\eta_{k,n}=\eta_{k,n-1}y_n,
\end{equation*}
and the Ore condition that
\begin{equation*}
\label{e63}
\eta_{j,n-1}=\eta_{k,n-1}.
\end{equation*}
Similar  argument implies in $n-k$ steps the equalities
\[
\eta_{j,n}=\eta_{k,n}, \;\; \eta_{j,n-1}=\eta_{k,n-1},\;\; \cdots,\;\;
\eta_{j,k+1}=\eta_{k,k+1}.
\]
Now the relations
\[
\zeta_{j,k+1}\eta_{j,k+1}=\eta_{j,k}x_{k+1},   \quad
\zeta_{k,k+1}\eta_{k,k+1}=x_kx_{k+1}
\]
and the Ore condition again imply  $\eta_{j,k}= x_k.$ By
(\ref{zetaj}) we have
\begin{equation*}
\label{e64}
(\zeta_{j,k}\eta_{j,k}=\eta_{j,k}x_k) \in \Re _0
\end{equation*}
This is impossible, since $\eta_{j,k}= x_k,$ and the relations in $\Re _0$ are square-free.

We have shown that the assumption $\eta_{j,n}=\eta_{k,n}$, for
some $k>j$, leads to a contradiction. This proves (\ref{r64})
for all $j, 1 \leq j\leq n-1$.

We set
\begin{equation*}
\label{eta6} \eta_1=\eta_{1,n},\quad \eta_2=\eta_{2,n}, \; \cdots, \;
\eta_{n-1}=\eta_{n-1,n}.
\end{equation*}

Next  we prove (\ref{r62}).

By the inductive assumption we have
\begin{equation*}
\label{e64} \zeta_{k,k+1}\zeta_{k,k+2}\cdots\zeta_{k,n} = x_{k+1}
\cdots x_n \in \Ncal _0
\end{equation*}
and
\begin{equation*}
\label{e65} x_{k+1} \cdots x_n.\eta_{k+1}=x_k\cdots x_n.
\end{equation*}
Applying the relations \ref{zetaj} we obtain
\begin{equation*}
(\zeta_{j,j+1}\zeta_{j,j+2}\cdots
\zeta_{j,n})\eta_j=x_jx_{j+1}\cdots x_n \in \Ncal _0.
\end{equation*}
Denote the normal form $Nor_{\Re
_0}(\zeta_{j,j+1}\zeta_{j,j+2}\cdots \zeta_{j,n})$ modulo the
Gr\"obner basis $\Re _0, $ as
\begin{equation*}
\label{e66} \upsilon _j=Nor_{\Re
_0}(\zeta_{j,j+1}\zeta_{j,j+2}\cdots \zeta_{j,n}),
\end{equation*}
clearly,  $\upsilon _j \in \Ncal _0.$

We have to show that there is an equality of words in $\langle X
\rangle$.
\begin{equation*}
\label{e666} \upsilon _j = x_{j+1}x_{j+2} \cdots x_n.
\end{equation*}
The equality,
\begin{equation*}
\label{e667}  \upsilon _j\eta_j,=x_jx_{j+1}\cdots x_n \in \Ncal
\end{equation*}
implies
\begin{equation*}
\label{e668}
Nor(\upsilon_j \eta_j)=x_jx_{j+1}\cdots x_n,
\end{equation*}
as words in the free semigroup $\langle X \rangle.$ Furthermore
$\upsilon _j$ does not contain subwords of the type $xx,$ (this
can be easily seen using the weak cyclic condition). Thus
\begin{equation*}
\label{e669} \upsilon _j = x_{j_1}x_{j_2}\cdots x_{j_{n-1}}, \
\text{where} j_1 < j_2 \cdots <j_{n-1} \leq n.
\end{equation*}
and therefore
\begin{equation}
\label{e670} j_1 \leq j+1 .
\end{equation}
The theory of Gr\"obner basis implies the following relations  in
$\langle X \rangle.$
\begin{equation*}
\label{e671} x_jx_{j+1}\cdots x_n = Nor(\upsilon _j \eta_j)\preceq
\upsilon _j \eta_j = x_{j_1}x_{j_2}\cdots x_{j_{n-1}}\eta_j,
\end{equation*}
therefore $j \leq j_1.$ By the last inequality, and (\ref{e670})
only two cases are possible: \textbf{a.} $j_1= j;$ and\textbf{ b.} $j_1=j+1.$ Assume
that $j_1=j.$ It follows then that
\begin{equation*}
\label{e672}  \upsilon _j = x_j\cdots x_{k-1} x_{k+1} \cdots x_n,
\end{equation*}
for some $k, k\geq j+1$ (In the case when $k=j+1,$ $\upsilon_j =
x_j x_{j+2} \cdots x_n).$ Thus the equalities
\begin{equation*}
\label{e673}  \upsilon _j \eta_k = x_j\cdots x_{k-1} (x_{k+1}
\cdots x_n\eta_k)=x_jx_{j+1}\cdots x_n =\upsilon_j \eta_j,
\end{equation*}
hold in $\Scal _0.$  So, by the cancellation low in $\Scal _0,$
we obtain
$\eta_k= \eta_j,$ with $j < k,$ which is impossible. Hence the
assumption $j_1=j$ leads to a contradiction. This verifies
$j_1=j+1,$ which implies $\upsilon _j = x_{j+1} \cdots x_n,$
and therefore the desired equality
\begin{equation*}
x_{j+1} \cdots x_n \eta _j= x_{j} \cdots x_n
\end{equation*}
holds in $\Scal _0.$ The lemma has been proved.
\end{proof}

\begin{proofproposition3}
%\ref{regular5}
Condition (\ref{regular51}) is obvious. Lemma \ref{regular6} proves
\ref{regular52}, \ref{regular53}. By the choice of $\eta _i, \; 1
\leq i \leq n-1,$ the following equalities hold in $\Scal _0:$
\begin{equation*}
\begin{array}{ll}
x_n\eta_{n-1}\eta_{n-2}\cdots \eta_j&= x_{n-1}x_n\eta_{n-2}\cdots
\eta_j \\
&= \cdots \cdots\cdots\cdots\cdots\\
&= x_{j+1}\cdots x_{n-1}x_n\eta_j\\
& = x_{j+1}\cdots
x_{n-1}x_n,
\end{array}
\end{equation*}
which implies (\ref{regular54}). The proof of conditions
(\ref{regular55}), (\ref{regular56}), and (\ref{regular57}) is analogous
to the proof of (\ref{regular52}), (\ref{regular53}), and
(\ref{regular54}), respectively. It follows from the weak cyclic
condition, that the normal form $Nor(u)$ of a monomial $u \in
\langle X \rangle$, with the shape $u=ayyb, y \in X$ has the shape
$Nor(u) = a_1xxb_1 \in \Ncal_0, x\in X.$ Therefore $W$ is the
principal monomial of $\Scal _0.$ Condition (\ref{regular58}) is
obvious.
\end{proofproposition3}
The following  lemma is used for the Frobenius property.
\begin{lemma}\label{mainlemma}
For any monomial $u \in \Ncal^{!}$ there exist uniquely determined
$u^{\prime}$ and $u^{\prime \prime}$ in $\Ncal^{!}$, such that
\begin{equation}
\label{pairing} u \ast u^{\prime} = x_1x_2\dots x_n, \quad u^{\prime
\prime} \ast u = x_1x_2\dots x_n.
\end{equation}
\end{lemma}
\begin{proof}
Let $u$ be an element of $\Ncal^{!}$. Then
\begin{equation*}
u = x_{1}^{\varepsilon_{1}}\cdots x_{n}^{\varepsilon_{n}}
\end{equation*}
where for all $i$, $1 \leq i \leq n,$ one has $0 \leq
\varepsilon_{i} \leq 1.$ Let  $\eta_i, \theta_j,\; 1 \leq i, j-1
\leq n-1,\;$ be as in  Proposition \ref{regular5}.
Let\begin{equation*}
\begin{array}{lll}
 u^{\prime}&=&x_n^{(1-\varepsilon_n)}\ast
\eta_{n-1}^{(1-\varepsilon_{n-1})} \ast \dots \ast
\eta_{1}^{(1-\varepsilon_1)}.\\
&&\\
u^{\prime\prime}&=&\theta_n^{(1-\varepsilon_n)}\ast
\theta_{n-1}^{(1-\varepsilon_{n-1})} \ast \dots \ast
\theta_{2}^{(1-\varepsilon_2)}x_1^{(1-\varepsilon_1)}.\\
&&
\end{array}
\end{equation*}
It is easy to verify that the equalities \ref{pairing} hold. The
uniqueness of $u^{\prime}$ and $u^{\prime \prime}$ follows from
the cancellation law in $\Scal _0.$
\end{proof}
\begin{prooftheoremA}
Let $\Fcal$ be the quadratic algebra from the hypothesis of
Theorem A.
%(\ref{frobeniustheoremF}).

Then Lemma \ref{frobeniuslemma1} and Corollary \ref{frobeniuscorrolary2}
imply conditions (1) and (2) of the theorem.  For $0 \leq i$ we set
\begin{equation}\label{frobeniuseq1}
\begin{split}
\Ncal^{!}_i &= \{u \in \Ncal^{!} \mid   u \;\text{has length}\; i \} \\
\Fcal_i &= Span_\textbf{k}\ \Ncal^{!}_i.
\end{split}
\end{equation}

It is clear that $\Fcal_0 = \textbf{k}$, $\Fcal_i = 0$, for $i >
n$ and for $1 \leq i \leq n$ one has
\[\begin{array}{c}
dim_\textbf{k}\ \Fcal_i = \sharp\ \Ncal^{!}_i = \binom{n}{i},\\
 \text{in \
particular}, \ dim_\textbf{k} \Fcal_n = 1 .
\end{array}
\]

Clearly, $\Fcal$ is
graded : $\Fcal= \bigoplus_{0\leq i \leq n}\Fcal_i$, $\Fcal_i= 0,
$ for $i> n$.

It follows from Lemma \ref{mainlemma} that the map
\[
(-,-): \Fcal_i
\times \Fcal_{n-i} \rightarrow \Fcal_n
\] defined by
\[(u,v): = \; \text{ the
normal form of}\;  uv \;\text{ in $\Fcal$}
\]
 is a perfect duality. This proves
Theorem A.
\end{prooftheoremA}
Now we can prove Theorem \ref{Frobeniusth}
\begin{prooftheorem}
Let $A$ be a binomial skew polynomial ring. By  Fact \ref{priddykoszul}
every algebra with quadratic Gr\"{o}bner basis is Koszul, this
implies the Koszulity of  $A$. Furthermore from \cite{Anick} one
deduces that for every graded $\textbf{k}$-algebra $\Bcal$ with
quadratic Gr\"{o}bner basis, Anick's resolution of $\textbf{k}$ as
a $\Bcal$-module is minimal. We shall use now the terminology of
Anick. The set of obstructions (i.e. the leading monomials of the
elements of the reduced Gr\"{o}bner basis) for a binomial skew
polynomial $A$  is $\{x_jx_i \mid 1 \leq i< j\leq n \}$. Therefore
the maximal $k$ for which there exist $k$-chains is $k=n-1.$ In
fact the only $n-1$-chain is $x_nx_{n-1}\cdots x_1.$ It follows
then from a theorem of Anick, \cite{Anick}, that $gl.dim A = n$.
We have shown that $A$ is a Koszul algebra of finite global
dimension. Clearly, $A$ has polynomial growth. Furthermore, by Theorem A, the Koszul dual $A^{!}$ is
Frobenius. It follows then from  \ref{factkoszul} that $A$ is
Gorenstein, and therefore $A$ is Artin-Schelter regular.
\end{prooftheorem}
\begin{prooftheoremB}
Let $A= \textbf{k}\langle X \rangle/(\Re)$ be a quantum binomial
algebra. The implication \ref{theoremB1} $\Longrightarrow$
\ref{theoremB2} follows from Theorem \ref{regular2}. Assume now
that $A$ is binomial skew polynomial ring. By remark
\ref{skewimpliescycliccondition} $A$ satisfies the cyclic
condition (see also \cite{T96}), and therefore it satisfies  the weak cyclic condition. By Corollary \ref{frobeniuscorrolary2} the Koszul dual
$A^{!}$ is Frobenius, and has regular socle. This proves the
implication \ref{theoremB2} $\Longrightarrow$ \ref{theoremB1}.

The equivalence of conditions (\ref{theoremB2}) and
(\ref{theoremB3}) follows from Theorem \ref{skewimpliesYBE} (see
also \cite{T06} Theorem 9.7).

We have shown that conditions (1), (2), and (3) are equivalent.

Now it is enough to show that every  binomial skew polynomial ring
$A$ satisfies the conditions (a) , ..., (e).  Conditions (a) and
(b)  are clear. We have shown that $A$ is Artin Schelter regular.
It is shown in \cite{TM}, Corollary 1.6 that $A$ is a domain. It
is proven in \cite{T96} (see also \cite{T1}, and \cite{TM}) that
$A$ is left and right Noetherian. It follows from \cite{TJO} that
$A$ satisfies polynomial identity. Now as a finitely generated PI
algebra, $A$ is catenary, see \cite{Sch}.
\end{prooftheoremB}

{\bf Acknowledgments}.  This paper combines new and some non
published  results which were found during my my visits at MIT
(1994-95) and at Harvard (2002). I express my gratitude to Mike
Artin, who inspired my research in this area, for his
encouragement and
moral support through the years. My cordial
thanks to Michel Van Den Bergh for our stimulating and productive
cooperation, for drawing my attention to the study of
set-theoretic solutions of the Yang-Baxter equation. It is my
pleasant duty to thank  David Kazhdan for inviting me to Harvard, for our valuable and
stimulating discussions and for his continuous support through the years.

%\bibliography{mybibs}
%\bibliographystyle{amsabbrv}
\ifx\undefined\bysame
\newcommand{\bysame}{\leavevmode\hbox to3em{\hrulefill}\,}
\fi

\end{document}